\input eplain 
\magnification  1200
\ifx\eplain\undefined \input eplain \fi

\baselineskip13pt

\overfullrule=0pt


\def\wtilde{\widetilde} 
\font\symb=cmsy5

\font\scaps=cmcsc10

\font\smalsmalbf=cmbx8
\font\tenbi=cmmib10

\font\eightbi=cmmib9

\font\fivebi=cmmib5
\newfam\bmifam\textfont\bmifam=\tenbi\scriptfont
\bmifam=\eightbi\scriptscriptfont\bmifam=\fivebi

\font\smalltenrm=cmr8
\font\smallteni=cmmi8
\font\smalltensy=cmsy8
\font\smallsevrm=cmr6   \font\smallfivrm=cmr5
\font\smallsevi=cmmi6   \font\smallfivi=cmmi5
\font\smallsevsy=cmsy6  \font\smallfivsy=cmsy5
\font\smallsl=cmsl8      \font\smallit=cmti8

\def\smallfonts{\lineadj{80}\textfont0=\smalltenrm  \scriptfont0=\smallsevrm
                \scriptscriptfont0=\smallfivrm
    \textfont1=\smallteni  \scriptfont1=\smallsevi
                \scriptscriptfont0=\smallfivi
     \textfont2=\smalltensy  \scriptfont2=\smallsevsy
                \scriptscriptfont2=\smallfivsy
      \let\it\smallit\let\sl\smallsl\smalltenrm}

\font\eightbi=cmmib10 at 8pt

\font\smathbold=msbm8\font\ssmathbold=msbm5
\font\mathbold=msbm9 at 10pt\def\Hn{{\hbox{\mathbold\char72}}^n}\def\Sn{{S^{n-1}}}
\def\sHn{{{\hbox{\smathbold\char72}}_{}^n}}\def\ssHn{{{\hbox{\ssmathbold\char72}}_{}^n}}
\def \C{{\hbox{\mathbold\char67}}}
\def\S{{\cal S} }
\def\R{{\hbox{\mathbold\char82}}}\def\Z{{\hbox{\mathbold\char90}}}
\def\sR{{\hbox{\smathbold\char82}}}

\def\avg{{-\hskip-1em\int}}

\def\e{\epsilon}

\def\P{{\cal P}}


\def\imp{{\;\;\Longrightarrow\;\;}}
\def\lineadj#1{\normalbaselines\multiply\lineskip#1\divide\lineskip100
\multiply\baselineskip#1\divide\baselineskip100
\multiply\lineskiplimit#1\divide\lineskiplimit100}

\def\remark#1.{\medskip{\noin\bf Remark #1.\enspace}}
\def\endpf{$$\eqno/\!/\!/$$}

\def\pf#1.{\smallskip\noin{\bf  #1.\enspace}}

\def\noin{\noindent}
\def\D{{\cal D}}

\def\ds{\displaystyle}
\def\ts{\textstyle}

\def\e{\epsilon}\def\part{\partial_t}

\def\isn{\int_{S^{n-1}}}\def\wtilde{\widetilde}
\def\what{\widehat}
\def\Rn{\R^n}

\def\irn{\int_{\sR^n}}
\def\ref#1{{\bf{[#1]}}}

\def\p{\partial}

\def\half{{\ts{1\over2}}}

\font\symb=cmsy9
\def\lfl{{\hbox{\symb\char 98}}}\def\rfl{{\hbox{\symb\char 99}}}

\def\s{\sigma}
\def\supp{{\rm supp }\,}
\def\H{{\hbox{\mathbold\char72}}}\def\b{\beta}
\def\a{\alpha}
\def\na{{n/\a}}\def\nsa{{n\over \a}}\def\an{{\a/n}}\def\asn{{\a\over n}}
\def\nna{{n\over n-\a}}
\def\bp{{\b'}}\bigskip\def\b{\beta}\def\bi{{1\over\b}}\def\bpi{{1\over\b'}}\def\ga{\gamma}
\def\nf{{\|f\|_\na^\na}}\def\ntf{{\|T_g f\|_\na^\na}}
\centerline{\bf  Sharp Adams and Moser-Trudinger inequalities on $\R^n$}
\centerline{\bf and other spaces of infinite measure}
\bigskip \centerline{Luigi Fontana, Carlo Morpurgo}
\footnote{}{\smallfonts
\hskip-2.4em  This work was partially supported by NSF Grant DMS-1401035 and by Simons Foundation Collaboration Grant 279735}
\midinsert
{\smalsmalbf Abstract. }{\smallfonts We derive  sharp Adams inequalities for the Riesz and other potentials of functions with arbitrary
compact support in $\sR^n$.  Up to now such results were only known for a class of functions  whose supports have uniformly bounded measure.

We obtain several sharp Moser-Trudinger inequalities for the critical Sobolev space $W^{\a,\nsa}$  on $\sR^n$ and on the hyperbolic
space $\sHn$. The  only known results so far are for $\alpha = 1$, both on $\sR^n$ and $\sHn$, and for $\alpha = 2$ on $\sR^n$. Other sharp inequalities are obtained  for general elliptic operators with constant coefficients  and for trace type Borel measures. We introduce critical potential spaces on $\sR^n$ on which our results can be extended to noninteger values of $\alpha$.

}

\endinsert\medskip

\centerline{\bf Table of contents}\bigskip
\item{1.} {\scaps   Introduction}\smallskip
\itemitem{-} {\sl Theorem 1: Inequalities  \`a la Ruf, for the  Riesz potential  and the higher order gradients  on $\Rn$}
\itemitem{-} {\sl Corollary 2: Inequalities  \`a la Adachi-Tanaka, for the  Riesz potential  and the higher order gradients  on $\Rn$}
\itemitem{-} {\sl Theorem 3: Inequalities on measure spaces with arbitrary measure}
\itemitem{-} {\sl Corollary 4:  Inequalities for integrable potentials on $\Rn$}
\medskip
\item{2.} {\scaps  Adams inequalities for general homogeneous Riesz-like potentials}\smallskip
\itemitem{-} {\sl Theorem 7:  Adams inequalities \`a la Ruf}
\itemitem{-} {\sl Corollary 8: Adams inequalities \`a la Adachi-Tanaka}
\itemitem{-} {\sl Proof of Theorem 7: overview and preliminary lemmas}
\itemitem{-}{\sl Proof of Theorem 7}
\itemitem{-} {\sl A simple proof of an Adams inequality \`a la  Adachi-Tanaka}
\itemitem{-} {\sl Proof of Corollary 8}
\itemitem{-} {\sl Proof of Theorem 1}
\medskip
\item{3.}  {\scaps Proof of Theorem 3}\smallskip

\itemitem{-} {\sl Proof of Corollary 4}
\itemitem{-} {\sl Proof of Corollary 6}\medskip
\item{4.}  {\scaps Further consequences of Theorem 3 and Theorem 7}\smallskip

\itemitem{-} {\sl Theorem 15: Inequalities for elliptic operators with constant coefficients}
\itemitem{-} {\sl Theorem 18: Inequalities in hyperbolic space}\smallskip
\item{5.} {\scaps Further results and extensions}\smallskip

\itemitem{-} {\sl Inequalities in critical potential spaces}
\itemitem{-} {\sl Inequalities for more general  Borel measures} 
\itemitem{-} {\sl  A sharp Trudinger  inequality on bounded domains without boundary conditions}
\medskip\item{6.}  {\scaps Appendix}
\smallskip

\bigskip
\bigskip
\centerline{\scaps 1. Introduction}\bigskip

It is well-known that if $\Omega$ is a bounded open set in $\Rn$ with smooth enough boundary, and $\alpha$ is an integer such that $0<\a<n$, then there exist constants $\ga,C>0$ such that
$$\int_\Omega e^{ \ga |u(x)|^{\nna}} dx\le C\eqdef{*1}$$
for every function $u\in W^{\a,\nsa}(\Omega)$ such that $\|u\|_{W^{\a,\nsa}}\le 1$. Here $W^{\a,p}(\Omega)$ denotes the usual Sobolev space of functions $u\in L^p(\Omega)$ having all derivatives up to order $\a$ also in $L^p(\Omega)$. More generally, such result holds also for noninteger $\alpha$, provided we interpret $W^{\a,p}(\Omega)$ as the space of  Bessel potentials of order $\a$, restricted to $\Omega$; for a proof see the papers by Strichartz [St1, St2], and the celebrated paper by Trudinger [Tr] for the case $\a=1$.

The problem of finding the sharp version of \eqref{*1}  goes typically  as follows: {\it given a subspace $B$ of $W^{\a,\nsa}(\Omega)$, endowed with a norm $\|\cdot\|_B$ equivalent to the usual  norm in $W^{\a,\nsa}(\Omega)$, what is the largest interval of $\gamma$'s for which \eqref{*1} is valid uniformly, as $\|u\|_B\le 1$?}
The largest such $\gamma$, if it exists, is called the {\it best constant} for \eqref{*1}, relative to the space $(B,\|\cdot\|_B)$. 

The most famous result in this direction is due to Adams [A] who found the best constant for the space $B=W_0^{\a,\nsa}(\Omega)$, endowed with the norm $\|\nabla^\a u\|_{\na}$, where $\a$ is an integer, and $\nabla^\a$ is the higher order gradient, defined as 
$$\nabla^\alpha=\cases{(-\Delta)^{\a\over2} & if $\alpha$ even\cr 
\nabla(-\Delta)^{\a-1\over2} & if $\a$ odd.\cr}\eqdef{gradient}$$
 Adams' result generalized an equally famous result due to Moser nearly two decades  earlier [Mo], who found the best constant in the case  $B=W_0^{1,n}(\Omega)$. Adams' method relies on a sharp exponential inequality for the Riesz potential   
$$I_\a*f (x)=\irn|x-y|^{\a-n}f(y)dy.$$
Using O'Neil's Lemma and other technically challeging 1-dimensional estimates, Adams proves that there is $C>0$ depending only on $|\Omega|$, such that for every  $f$ with $\|f\|_\na\le 1$
$$\int_\Omega \exp\bigg[{1\over|B_1|}|I_\a* f(x)|^{\nna}\bigg] dx\le C,\qquad |B_1|={\omega_{n-1}\over n}\eqdef{*2}$$
where $\omega_{n-1}$ is the volume of the $(n-1)-$dimensional sphere, and where uniformity in $f$ is lost if $|B_1|^{-1}$ is replaced by a larger constant.
Using that  each  $u\in C_0^\infty(\Omega) $ (or each $u\in W_0^{\a,\nsa}(\Omega)$) can be  estimated sharply in terms of  the Riesz potential of $|\nabla^\a u|\in L^{\nsa}(\Omega)$, via the identities
$$u=\cases{c_\a I_\a*(-\Delta )^{\a\over2}u, &  if $\alpha$  even \cr 
c_{\a+1} \nabla I_{\a+1}*\nabla(-\Delta )^{{\a-1}\over2}u &  if $\alpha$  odd \cr}
\qquad c_\a={\Gamma\big({n- \a\over2}\big)\over 2^\a\pi^{n/2}\Gamma\big({ \a\over2}\big)}\eqdef{identities}$$ 
Adams finds that the best constant  for \eqref{*1}, when  $B=W_0^{\a,\nsa}(\Omega)$ is given by 
$$\gamma_{n,\a}:=\cases{\ds{c_\a^{-\nna}\over|B_1|}, &  if $\alpha$  even  \cr\ds{\big((n-\a-1)c_{\a+1}\big)^{-{\nna}}\over|B_1|} & if $\alpha$ odd.\cr}\eqdef{gamma}$$
 
Inequalities such as \eqref{*1} and  \eqref{*2} are  referred to as Moser-Trudinger (MT) and Adams  inequalities, respectively.

Since Adams' work, countless papers have been published on sharp  Moser-Trudinger and Adams inequalities in various forms and settings. Often, such results were motivated by questions in conformal geometry or nonlinear PDE, but it is fair to say that over the years the challenge of finding best constants in inequalities such as \eqref{*1}  and \eqref{*2} has taken a life of its own, generating a rather active area of research. In particular, the sharp Moser-Trudinger inequality in the form \eqref{*1} is now very well understood on $W_0^{\a,\nsa}(\Omega)$ when $\Omega$ is  an  open subset with {\it finite measure} of a Riemannian or subRiemannian manifold. Likewise, Adams' inequality in the form (2) is completely understood even on arbitrary measure spaces with finite measure, where $I_\alpha * f$ is replaced by a general integral operator [FM1].


On the contrary, when $\Omega$ has {\it infinite measure}, relatively  few sharp Moser-Trudinger inequalities exist for the space $W_0^{\alpha,\nsa}(\Omega)$. More importantly,  {\it no version of the Adams inequality for the Riesz potential exists} even when $\Omega=\Rn$. In our opinion the literature dealing with Moser-Trudinger inequalities on $\Rn$ is somewhat confusing; below we will try to summarize what we consider the most important achievements on this subject. 

First off, if $|\Omega|=\infty$, inequality  \eqref{*1} is obviously false, and the exponential  needs to be  suitably regularized when $u$ is near $0$. The standard way to do this is to consider   
$$\exp_N(t)=e^t-\sum_{k=0}^N{t^k\over k!},\qquad N=0,1,...$$
and recast \eqref{*1} into the inequality
$$\int_\Omega \exp_{[\nsa-2]}\Big(\gamma|u(x)|^{\nna}\Big)dx\le C\eqdef{*3}$$
where $[x]$ denotes the ceiling of $x$, i.e. the  smallest integer greater or equal  $x$, for $x\in \R$. If $\Omega$ has infinite inradius it is not hard to find a family of functions $u$ with  $\|\nabla^\alpha u\|_\na=1$,  along which \eqref{*3} is false for any $\gamma>0$.  For such $\Omega$'s in general   $\|\nabla^\alpha u\|_\na$ is not a norm in $W_0^{\a,\nsa}(\Omega)$, however, using a Sobolev norm, inequality \eqref{*3} without sharp constants, was derived first by Ogawa [Og] when $n=2$, $\alpha=1$,  and by Ozawa [Oz] for arbitrary $n$ and $\alpha<n$, in the case $\Omega=\Rn$. In particular,  such results were obtained under the condition
$$\max\big\{\|u\|_{n},\,\|\nabla u\|_n\big\}\le 1,\qquad u\in W^{1,n}(\Rn)\eqdef{*4}$$
in [Og] and 
$$\max\big\{\|u\|_{\na},\,\|(-\Delta)^{{\a\over2}} u\|_\na\big\}\le 1,\qquad u\in W^{\a,\nsa}(\Rn)\eqdef{*5}$$
in [Oz]. Clearly the quantities on the LHS of \eqref{*4} and \eqref{*5} are norms in the corresponding Sobolev spaces, and can be replaced by any other equivalent norms, since the exponential constants are not sharp. Also, note that the results in [Og], [Oz] are valid on $W_0^{\a,\nsa}(\Omega)$ any open $\Omega\subseteq \R^n$, or on $W^{\a,\nsa}(\Omega)$, if $\Omega$ is smooth enough.

The sharp version of \eqref{*3} was first obtained for $\alpha=1$ by Cao [Cao] in dimension 2 and by Panda [Pa] and Do \' O [Do\'O] in arbitrary dimension, who proved that under condition  \eqref{*4}
$$\int_{\sR^n} \exp_{n-2}\Big(\gamma|u(x)|^{n\over n-1}\Big)dx\le C\eqdef{*6}$$
when $\gamma<\gamma_{n,1}=n\omega_{n-1}^{1\over n-1}$, the original Moser sharp constant. A couple of years later Adachi and Tanaka reproved Cao and Panda's inequality and cast it in the equivalent dilation invariant form 
 $$\int_{\sR^n} \exp_{n-2}\bigg[\gamma\bigg({|u(x)|\over \|\nabla u\|_n}\bigg)^{n\over n-1}\bigg]dx\le C\bigg({\|u\|_n\over\|\nabla u\|_n}\bigg)^n,\qquad u\in W^{1,n}(\Rn)\setminus\{0\}.\eqdef{*7}$$
 Moreover, Adachi and Tanaka [AT] verified  that the above inequality \eqref{*7} fails if $\gamma=\gamma_{n,1}$ by means of the usual Moser sequence, along which the left hand side is bounded away from 0 (in fact bounded above too by Moser's inequality), whereas the  right hand side tends to~0. From now on we will refer to \eqref{*6} under condition \eqref{*4}, or equivalently \eqref{*7},  as the Adachi-Tanaka inequality.

In  2008 Li-Ruf proved that \eqref{*6} holds at the critical constant $\gamma=\gamma_{n,1}$  for the usual  Sobolev norm, namely under the {\it Ruf condition}
$$\|u\|_n^n+\|\nabla u\|_n^n\le 1\eqdef{*8}$$
and verified that the inequality fails if $\gamma>\gamma_{n,1}$ along the usual Moser sequence [LR]. This result was first proved for $n=2$ by Ruf, three years earlier [Ruf]. With a simple dilation argument it is easily seen that Ruf's inequality implies 
Adachi-Tanaka's inequality (See proofs of Corollaries 2 and 8). In the context of the Heisenberg group the Li-Ruf result was obtained recently by Lam and Lu [LL1], and on general noncompact Riemannian manifolds  Yunyan Yang [Y] proved that the corresponding Moser-Trudinger inequality holds for  for $\gamma<\gamma_{n,1}$, under the Ruf condition \eqref{*8}.

In [LL2, Thm 1.5] Lam-Lu established \eqref{*6} for $\gamma=\gamma_{n,2}$ under the condition
$$\|u\|_{n/2}^{n/2}+\|\Delta u\|_{n/2}^{n/2}\le 1,\eqdef{LL}$$
which is the second order version of the Li-Ruf result. Lam-Lu obtained this result by reducing it to the sharp Moser-Trudinger inequality on bounded domains with homogeneous Navier boundary conditions, derived in [Tar].

 Recently, Ibrahim-Masmoudi-Nakanishi [IMN]  ($n=2$) and Masmoudi-Sani [MS1]  ($n\ge2$) derived the following more refined  inequality 
$$\int_{\sR^n}{\exp_{n-2}\big(\gamma_{n,1}|u(x)|^{n\over n-1}\big)\over(1+|u|)^{n\over n-1}}dx\le C,\qquad \|\nabla u\|_n\le 1,\;\;\|u\|_n\le 1\eqdef{IMN}$$
which incorporates  Li-Ruf's result. In [MS2], Masmoudi-Sani proved that for $\alpha=2,n=~4$
$$\int_{\sR^4}{\exp_{0}\big(\gamma_{4,2}^{} |u(x)|^2\big)\over(1+|u|)^{2}}dx\le C,\qquad \|\Delta u\|_2\le 1,\;\;\|u\|_2\le 1\eqdef{MS}$$
which implies (in dimension 4) the Lam-Lu sharp inequality  under condition \eqref{LL}. Even more recently, \eqref{MS} was extended to all dimensions  in [LTZ].

The higher order versions of Li-Ruf, Lam-Lu, or Masmoudi et al. results for $\alpha>2$ are not known, if the $W^{\alpha,\nsa}$ norm of $u$ is to be controlled only by the $L^{\nsa}$ norms of $u$ and $\nabla^\alpha u$. In particular, for integer $\alpha\in(0,n)$,   the best constant for the inequality
$$\int_{\sR^n} \exp_{[\nsa-2]}\Big(\gamma |u(x)|^{\nna}\Big)dx\le C\eqdef{*11}$$
under the Ruf condition
$$ \|u\|_{\na}^\na+\|\nabla^\a u\|_\na^\na\le 1,\eqdef{*12}$$
is not known, except for the cases discussed above,  $\alpha=1$, and $\alpha=2$.

More generally, the sharp Adams inequality for the Riesz potential on $\R^n$
$$\int_{\sR^n} \exp_{[\nsa-2]}\Big({1\over|B_1|}|I_\a*f(x)|^{\nna}\Big)dx\le C\eqdef{*13}$$
for all $f\in L^\nsa(\Rn)$ such that $I_\a*f$ is well defined and in $L^{\nsa}(\Rn)$ with  
$$ \|f\|_{\na}^\na+\|I_\a*f\|_\na^\na\le 1,\eqdef{*14}$$
 it is not known for {\it any} value of $\alpha\in (0,n)$. Arguing like in [A] using \eqref{identities}, it is clear that is $\alpha$ is even, \eqref{*13} and \eqref{*14} imply \eqref{*11} under condition \eqref{*12} with $\gamma=\gamma_{n,\a}$. 
The  relation between  \eqref{*14} and \eqref{*12}  is not so  trivial in the case $\alpha$ odd (see proof of Theorem 1). 

It is worth noting that the known results on sharp Moser-Trudinger inequalities under the Ruf condition or the Masmoudi-Sani type condition, ultimately rely on symmetrization tools such as  the P\'olya-Szeg\"o inequality and Talenti's comparison theorem, which are not suitable to derive the corresponding Adams inequalities, at least for arbitrary orders.

The Ruf norm in \eqref{*12} is in some sense minimal, in regard to  the number of derivatives that it  involves. 
 Under norm conditions more restrictive than \eqref{*12}   sharp higher order results do exist. Indeed, the sharp inequality in \eqref{*11} holds  for any $\alpha\in(0,n)$ under the condition  
$$\|(I-\Delta)^{\a\over2}u\|_\na\le 1.\eqdef{*15}$$
This result goes back to Adams, who proved it for $\alpha=2$  in his original 1988 paper [A, Thm. 3]. Strictly speaking, Adams proved that if  $\Omega$ an open and bounded set with measure $|\Omega|\le 1$,  $\alpha=2$ and if $u\in W^{2,{n\over2}}(\R^n)$ satisfies \eqref{*15}, then the basic Moser-Trudinger inequality \eqref{*1} holds at the critical index $\gamma_{n,2}$. To prove this result Adams  modified slightly the proof that he gave of \eqref{*2} for Riesz potentials on   $L^\nsa(\Omega)$, adapting it   to Bessel potentials on $L^{\nsa}(\R^n)$; in particular, he proved that for $\alpha=2$ there exists $C>0$, independent of $\Omega$, such that 
$$\int_\Omega \exp\bigg[{1\over|B_1]} |G_\a* f(x)|^{\nna}\bigg] dx\le C,\qquad \|f\|_{\na}\le 1\eqdef{*16}$$
where $G_\alpha$ denotes the usual Bessel potential. In this setting it is not necessary to have $\Omega$ be a bounded open set, it's enough that $\Omega$ be measurable with measure $\le 1$. Adams' proof of \eqref{*16} for $\alpha=2$, however,  is actually working for any $\alpha\in (0,n)$, after trivial modifications; it appears that the only reason why Adams considered $\alpha=2$ is because in this case the norm in \eqref{*15} coincides precisely with the usual full Sobolev norm on $W^{2,{n\over2}}$. The fact that the Bessel potential kernel $G_\a$ decays quite well  at infinity is what, 
at the end of the day, makes Adams' proof go through with few modifications from the finite measure case. We shall return to this point shortly.

It must be noted that  \eqref{*16}, valid for any measurable $\Omega$ with $|\Omega|\le1$,  implies that (and it is in fact equivalent to) 
$$\int_{\sR^n} \exp_{[\nsa-2]}\Big({1\over|B_1|}|G_\a*f|^{\nna}\Big)dx\le C,\qquad \|f\|_{\na}\le 1.\eqdef{*17}$$
The proof of this fact is rather straightforward: write $\R^n=\Omega\cup\Omega^c$ with $\Omega=\{x: |G_\a*f(x)|\ge1 \}$, and split the integral accordingly. Clearly, $|\Omega|\le 1$,  (since $\|G_\a*f\|_\na\le \|f\|_\na\le1$), and the integral over $\Omega^c$ can be estimated by $\|G_\a*f\|_\na$ by writing the exponential as a  Taylor series. This observation applies also to \eqref{*13}, \eqref{*11}, and all the other inequalities on $\Rn$ or on spaces with infinite measure. In short, the regularized exponential is nothing more than a gimmick, and can be replaced by the usual exponential provided that  the resulting  inequality holds on measurable sets of measure no greater than 1. In section 2  we will state an elementary ``Exponential Regularization Lemma", which will be used implicitly throughout this paper, in order to pass from inequalities over sets of finite measure to inequalities over the whole space.

Unaware of Adams' result, Ruf-Sani [RS] proved \eqref{*11} under the condition \eqref{*14} when $\alpha$ is an even integer,  not using the Bessel potential approach, but rather comparison theorems for  suitable Navier boundary value problems.
Recently, Lam-Lu [LL2]  presented a complete proof of the same result for any $\alpha\in (0,n)$, using the Bessel potential approach.


In this paper we derive both the sharp Adams and Moser-Trudinger inequalities on $\R^n$ for the full range of $\alpha$, under the Ruf conditions \eqref{*14} and \eqref{*12} respectively:

\proclaim Theorem 1. If  $0<\alpha<n$ there exists a constant $C=C(\alpha,n)$ such that:
 \smallskip\item{(a)} For every measurable and compactly supported  $f:\Rn\to\R$ such that 
$$\|f\|_\na^\na+\|I_\a*f\|_\na^\na\le 1\eqdef{100}$$
and for all measurable $E\subseteq \R^n$ with $|E|<\infty$, we have 
$$\int_E \exp\bigg[{1\over|B_1|}|I_\a*f(x)|^\nna\bigg]dx\le C(1+|E|)\eqdef{101}$$
and 
$$\int_{\sR^n}\exp_{[\nsa-2]}\bigg[{1\over|B_1|}|I_\a*f(x)|^\nna\bigg]dx\le C.\eqdef{101b}$$
\smallskip\item{(b)} If $\a$ is an   integer  then, for every  $u\in W^{\a,\nsa}(\Rn)$ such that 
$$\|u\|_\na^\na+\|\nabla^\a u\|_\na^\na\le1\eqdef{12b}$$
and for all measurable $E\subseteq \Rn$ with $|E|<\infty$, we have
$$\int_E\exp\Big[\gamma_{n,\a}|u(x)|^{n\over n-\a}\Big]dx\le C(1+|E|),\eqdef {12c}$$
and 
$$\int_{\sR^n}\exp_{[\nsa-2]}\Big[\gamma_{n,\a}|u(x)|^{n\over n-\a}\Big]dx\le C.\eqdef{12cc}$$
\item{(c)} If $|E|>0$, the exponential constants  in \eqref{101}, \eqref{101b}, \eqref{12c} and  \eqref{12cc} are sharp, that is, they cannot be replaced by  larger constants under the corresponding conditions \eqref{100}, \eqref{12b}.

\par

\def\U{{\cal U}}

As an immediate consequence of Theorem 1 we have  a version of the Adachi-Tanaka result, with sharp control on the right hand side constant. We will in fact consider a more general result for the family of Sobolev norms
$$\Big(\|u\|_{\na}^{q\na}+\|\nabla^\a u\|_{\na}^{q\na}\Big)^{\a\over qn},\qquad 1<q\le+\infty$$
for which the exponential inequality at the critical index fails.

\proclaim Corollary 2. If $\a$ is an integer, $0<\alpha<n$,  there exists $C=C(n,\a)>0$ such that for  $0<\theta<1$ and $1<q\le+\infty$, and  for all  measurable  $E\subseteq \Rn$ with $|E|<\infty$ 
$$\int_E \exp\Big[\theta\gamma_{n,\a}|u(x)|^\nna\Big]dx\le C( 1-\theta)^{-{1\over q'}}(1+|E|),\eqdef{102}$$
and
$$\int_{\sR^n}\exp_{[\nsa-2]}\Big[\theta\gamma_{n,\a}|u(x)|^\nna\Big]dx\le C( 1-\theta)^{-{1\over q'}}\eqdef{102b}$$
for all $u\in W^{\a,\nsa}(\Rn)$ such that 
$$\cases{\|u\|_\na^{q\na}+\|\nabla^\a u \|_\na^{q\na}\le 1, &  if $q<\infty$\cr\cr   \max\big\{\|u\|_\na,\,\|\nabla^\a  u\|_\na\big\}\le 1, & if  $q=+\infty.$\cr}\eqdef{103}$$
The inequalities \eqref{102},\eqref{102b} are sharp, in the sense that for given $q\in(1,\infty]$,  there exists a family of functions $u_\theta\in W^{\a,\nsa}(\Rn)$ satisfying \eqref{103} for which  \eqref{102} and \eqref{102b} are reversed; in particular, the exponential integral cannot be uniformly bounded if $\theta=1$.  
\par
See also  the recent  papers [CST] and  [LLZ], where the relations between the Li-Ruf and the Adachi-Tanaka type results as in \eqref{103}  are explored in more detail. 

A relevant aspect illustrated by  Corollary 2 is that the sharp inequalities \eqref{12c} and \eqref{12cc} of Theorem 1 do not  hold under \eqref{103} when $q>1$, whereas they do hold when $q\le 1$, since  the $q-$norms  in \eqref{103} are greater than the Ruf norm when $q<1$. This partly justifies the claim that among all possible natural Sobolev norms  in $W^{\a,\nsa}(\Rn)$ that involve  only $u$ and $\nabla^\a u$,  the Ruf norm  yields the least restrictive condition under which a sharp Moser-Trudinger inequality of the type \eqref{*11} holds at the critical constant $\gamma=\gamma_{n,\a}$.

Note that in the case $q=\infty$ one can immediately obtain a general dilation invariant version of Adachi-Tanaka analogous to \eqref{*7}:
$$\int_{\sR^n}\exp_{[\nsa-2]}\bigg[\theta\gamma_{n,\a}\bigg({|u(x)|\over\|\nabla^\a u\|_\na}\bigg)^\nna\Big]dx\le {C\over( 1-\theta)}\bigg({\|u\|_\na\over\|\nabla^\a u\|_\na}\bigg)^\nsa,$$
valid for all $u\in W^{\a,\nsa}(\R^n)\setminus\{0\}$. This is achieved by  first by replacing $u$ with $u/\|\nabla^\a u\|_\na$ in \eqref{102b},  and then by making the change of variable $x\to\lambda x$, with $\lambda^\a=\|u\|_\na/\|\nabla^\a u\|_\na$.\smallskip

Another perhaps interesting observation, is that Corollary 2 implies the following weak Masmoudi-Sani type inequalites:

$$\int_{\sR^n}{\exp_{[\nsa-2]}\Big[\gamma_{n,\a}|u(x)|^\nna\Big]\over(1+|u(x)|)^{\nna}}dx\le Cq, \eqdef{masmoudi}$$
valid for $1\le q<\infty$ and for all $u$ such that $\|u\|_\na^{qn/\a}+ \|\nabla^\a u\|_\na^{q\na}\le 1$, and
$$\int_{\sR^n}{\exp_{[\nsa-2]}\Big[\gamma_{n,\a}|u(x)|^\nna\Big]\over(1+|u(x)|)^{\nna(1+\e)}}dx\le {C\over\e}, \eqdef{masmoudi2}$$
for all $\e>0$, and for  $u$ such that $\|u\|_\na\le 1$ and $\|\nabla^\a u\|_\na\le 1$. These inequalities can be obtained by taking $E=\{u\ge 1\}$ in \eqref{102}, and integrating in $\theta$ directly, or after multiplying by  $(1-\theta)^\e$. Better inequalities can be obtained with this method, but still without reaching the optimal Masmoudi-Sani type of result.

A few remarks are in order here regarding the Riesz potential as an operator on the critical $L^p$ space. It is easy to see that for $f$  compactly supported and in $L^\nsa(\R^n)$ the potential  $I_\a*f$ is well defined almost everywhere, and for large $x$ it is a continuous function  that  decays at least as $|x|^{\a-n}$. However, the Riesz potential cannot be extended to a continuous operator on $L^\nsa$. In fact, consider the sequence of functions $f_k(x)=\chi_{2\le|x|\le k}^{}(x)|x|^{-\a}/\log|x|$,   which converges  in $L^\nsa$ to a function $f$, but pointwise $I_\a*f_k(x)\to~+\infty$  (hence $I_\a*f\equiv +\infty $ on $\R^n$.) It is however possible to extend the class of functions for which a potential is defined, by first adding the requirement that such potentials  be in $L^{\nsa}$, and then by constructing the smallest closed extension of the Riesz potential, as a densely defined operator from $L^\nsa$ to itself. Details of the construction of these {\it critical potential spaces} will be given in section 5, where we state a version of Theorem 1 valid  for  all $\alpha\in(0,n)$. For $\alpha<n/2$ these spaces  coincide with the  usual Bessel potential spaces, but it is unclear whether or not this  holds for all $\alpha<n$.

Theorem 1 was our main motivation, however in this paper we will obtain more general results at a minimal extra cost. Specifically, we will obtain a version of Theorem~1 for convolution operators of form 
$$T_gf(x)=\int_{\sR^n}g(x-y) f(y) dy\eqdef{104}$$
where $g(x)$ is a {\it nonnegative homogeneous function of order $\a-n$}, in the same spirit as the results in [FM1] and [CL] in the context  of sets of finite measure. The sharp constant in this case is given by 
$${n\over{\isn|g(\omega)|^{\nna}d\omega}}\eqdef{105}$$
provided that $g$ is smooth enough. This result will be also derived for vector-valued kernels $g(x-y)$, whose components are nonnegative and homogeneous of order $\a-n$; in this case the operator in \eqref{104} will be acting on vector-valued functions, and the sharp constant is still given as in \eqref{105} (see Theorem 7).
We will use this more general Adams inequality to obtain the Moser-Trudinger inequality in Thm 1 b) for general  invertible elliptic operators which are {\it homogeneous} of order $\a<n$  (see Theorem 15).

Our method also allows us to obtain, with minimal modifications,  the results in the form \eqref{101}, \eqref{12c} when the non-regularized exponential is integrated against a general Borel measure $\nu$ satisfying
$$\nu\big(B(x,r)\big)\le Q r^{\sigma n},\qquad \forall x\in \Rn,\;\;r>0$$
for some $Q>0$ and $\sigma\in (0,1]$.
The resulting Adams/Moser-Trudinger inequalities analogous to \eqref{101} and \eqref{12c} are of trace type, and the sharp constant in this case is $\sigma \gamma_{n,\a}$. Similar results hold for the general operators as in \eqref{104}. Examples of such measures are the Hausdorff measures on  submanifolds of  $\Rn$. Another example  is the ``singular measure" with density $|x|^{(\sigma-1)n}$, which was considered for example in [LL2] and [AY], and in other papers dealing with domains of finite measure.  For those measures, and even for more general ones,  we can even obtain the results in the full regularized form over the whole $\Rn$, ans in \eqref{101b}, \eqref{12cc}. For sake of clarity we will not state our main results for such more general measure, however  in section 5  we will provide a few  explicit statements.

In this paper we also treat the case of non-homogeneous potentials and operators;  ironically, these are considerably easier to handle, given the current abstract Adams' machinery developed in [FM1]. In an effort to understand the results in  [RS] for the Bessel operator $(I-\Delta)^{\a\over2}$, and unaware of the results in [LL2], we soon realized that the abstract results in our paper [FM1] could have been easily extended to treat spaces of infinite measure, by simply adding one integrability hypothesis.
Specifically, suppose that  $(M,\mu)$ and $(N,\nu)$ are measure spaces  and suppose that  that  $T$ is an integral operator of type
$$Tf(x)=\int_M k(x,y)f(y)d\mu(y),\qquad x\in N,$$
with $k:M\times N\to[-\infty,\infty]$ a $\nu\times\mu$-measurable function. Define
$$k_1^*(t):=\sup_{x\in N} [k(x,\cdot)]^*(t),\qquad
k_2^*(t):=\sup_{y\in M} [k(\cdot,y)]^*(t)$$
where $k[(x,\cdot)]^*$ is the nonincreasing rearrangement of $k(x,y)$ with respect to $y$ for fixed $x$, and $k[(\cdot,y)]^*$ is its analogue for fixed $y$. Then the following theorem holds:

\proclaim Theorem {3} (Refined Adams inequality on measure spaces). Let  $\nu(N)<\infty$ and  
$$k_1^*(t)\le A^{1\over\beta}t^{-{1\over\beta}}\big(1+H(1+|\log t|)^{-\gamma}\big),\qquad 0<t\le1\eqdef{106}$$
$$k_2^*(t)\le Bt^{-{1\over\sigma \beta}},\qquad t>0\eqdef{107}$$
$$J:={1\over A}\int_{1}^\infty\big (k_1^*(t)\big)^{\b}dt<\infty,\eqdef{108}$$
for some $\b>1$ and $0<\sigma\le 1$.
Then,  $Tf$ is finite a.e. for $f\in L^\bp(M)$, and  there exists a constant $C=C(\b,\sigma,\gamma,A,B,H)$ such that for each $f\in L^{\b'}(M)$ with $\|f\|_{\b'}\le 1$
$$\int_N \exp\bigg[{\sigma\over A}|Tf(x)|^\beta\bigg]d\nu(x)\le  Ce^{\sigma J}\big(1+J+\nu(N)\big).\eqdef{109}$$
If the condition $\nu(N)<\infty$ is dropped, then for all $f\in L^\bp(M)$ such that $\|f\|_\bp\le1$  we have  
$$\int_N \exp_{[\b'-2]}\bigg[{\sigma\over A}|Tf(x)|^\beta\bigg]d\nu(x)\le  Ce^{\sigma J}\big(1+J+\|Tf\|_\bp^\bp\big).\eqdef{110}$$
\par

The proof of this theorem amounts to some straightforward modifications in the proof of Theorem~1 in [FM1], using the new added  condition \eqref{108}, which effectively replaces $\mu(M)<\infty$. We also decided to track down the the dependance of the right-hand side in terms of $J$ and $\nu(N)$, which is used later in Corollaries 4 and 6, and might turn out to be useful in other instances.

The most obvious application of Theorem 3 is to convolution operators on $\Rn$:

\proclaim Corollary 4. Let $0<\a<n$ and suppose that $K_\a:\Rn\setminus\{0\}\to\R$ satisfies the conditions
$$K_\a(x)=g(x^*)|x|^{\a-n}+ O(|x|^{\a-n+\e}),\qquad |x|\le R,\; x^*={x\over|x|}\eqdef{111}$$
$$K_\a\in L^{\nna}\cap L^\infty\big(|x|\ge R\big)\eqdef{112}$$
for some $\e,R>0$, where $g\in L^{\nna}(S^{n-1})$. Then $K_\a*f$ is finite a.e for $f\in L^{\nsa}$, and   there exists $C>0$ such that for all $f\in L^{\nsa}(\Rn)$ with $\|f\|_\na\le 1$, and for each measurable $E\subseteq \Rn$ with $|E|<\infty$
$$\int_E\exp\bigg[{1\over A_g}|K_\a*f(x)|^\nna\bigg]dx\le  C(1+|E|),\eqdef{112a}$$
where $$A_g={{1\over n}{\isn|g(\omega)|^{\nna}d\omega}}.\eqdef{114}$$
Moreover, 
$$\irn\exp_{[\nsa-2]}\bigg[{1\over A_g}|K_\a*f(x)|^\nna\bigg]dx\le  C\big(1+\|K_\a *f\|_\na^\na\big).\eqdef{113}$$
for all $f\in L^\nsa(\R^n)$ such that $\|f\|_\na\le 1$.
If $g$ is smooth, then the exponential constant in \eqref{112a} (if $|E|>0$) and in \eqref{113} is sharp.\par

Note that the ``big O" notation  in \eqref{111} means that  $|O(|x|^{\a-n+\e})|\le C|x|^{\a-n+\e}$, for $|x|\le R$.

The proof of the inequality of the above corollary requires only to verify that conditions \eqref{111} and \eqref{112}  imply \eqref{106}, \eqref{107} and \eqref{108}, with $k(x,y)=K_\a(x-y)$, and $\beta=\nna$. Some of this has already been done in [FM1] in the context of finite measure spaces ([FM1, Theorem 8, Lemma 9]).
 The crucial point here, however,  is the integrability condition in \eqref{112}, which appears to be the bare required  minimum in order for the Adams' machinery to work in the infinite measure case, in addition to the usual local behavior of the the kernel. Condition \eqref{112} is what also markedly singles out the purely homogeneous case $K_\a(x)=g(x^*)|x|^{\a-n}$, treated separately in Theorem 7. 
 
Generally speaking  non-homogeneous invertible elliptic  operators will have a kernel satisfying \eqref{112}, and for those operators a sharp Moser-Trudinger inequality will hold.  As a first example, consider the Bessel potential $(I-\Delta)^{\a\over2}$, whose fundamental solution $G_\a(x)$ behaves like the Riesz potential locally, and decays exponentially at infinity. In fact an immediate consequence of Corollary 4, and the fact that $\|G_\a* f\|_p\le \|f\|_p$, is the aforementioned results by [A, Thm. 3],  [LL2] and [RS]:

\eject 

\proclaim Corollary 5. If $0<\a<n$ then there exists $C$ such that for all $u\in H^{\a,\nsa}(\Rn)$ so that 
$$\|(I-\Delta)^{\a\over2}u\|_\nsa\le 1$$
we have
$$\irn\exp_{[\nsa-2]}\bigg[{c_\a^{-\nna}\over|B_1|}|u(x)|^\nna\bigg]dx\le  C.\eqdef{bessel}$$
and the exponential constant is sharp.\par

In this paper we  define $H^{\a,\nsa}(\Rn)$ to be the space of Bessel potentials, which is known to coincide with the standard Sobolev space for $\alpha$ integer. In section 5, Theorem 21, we will give a version of Corollary 5  where $u$   is subject to  the Ruf condition $\|u\|_\na^\na+\|(-\Delta)^{\a\over2}u\|_\na^\na\le 1$.

In section 4, Theorem 15 we will show that if $P$  is any  non-homogeneous, elliptic, invertible, linear partial differential  operator with constant coefficients,  under the assumption that  $\a\ge n/2$ its  inverse has a  kernel satisfying \eqref{111} and \eqref{112}, and therefore a sharp Moser-Trudinger inequality holds for such $P$.

Yet another important application of Theorem 3 is in the hyperbolic setting: in section 4 we will derive a sharp Moser-Trudinger inequality in the full hyperbolic space for the higher order gradients,  extending recent results by Mancini-Sandeep-Tintarev [MST] and Lu-Tang [LT1], [LT2].

It should be remarked here that Corollary 4 can be likely  formulated so as to accommodate  more general (non-convolution) kernels satisfying
$$K_\a(x,y)=g(x,(x-y)^*)|x-y|^{\a-n}+O(|x-y|^{\a-n+\e}),$$
together with suitable integrability and boundedness conditions at infinity, in the same spirit as in [FM1, Thm. 8].

Another easy consequence of Theorem 3 in the case of homogeneous Riesz-like potentials on finite measure spaces, is the following slightly more general  version of the sharp Adams trace inequality given in [FM1]:

\proclaim Corollary 6.  a) Let  $\nu$ be a positive Borel measure on $\Rn$. Assume that  there exist  $Q>0$ and $\sigma\in (0,1]$ such that 
 $$\nu\big(B(x,r)\big)\le Q r^{\sigma n},\qquad \forall x\in \Rn,\,\forall r>0.\eqdef{b}$$
	 Let $g$ be homogeneous of order $\a-n$ on $\R^n$   ($0<\a<n$)  and with $g_{/{S^{n-1}}}\in L^{\nna}(S^{n-1})$. \smallskip
If $T_g f=g*f$, then there exists  $C=C(n,\a,\s,Q)$ such that for all $E,F\subseteq \Rn$ with $\nu(E)<\infty,\, |F|<\infty$ and for all $f\in L^{\na}(\Rn)$ with $\|f\|_\na\le 1$ and $\supp f\subseteq F$ we have
$$\int_E \exp\bigg[{\sigma\over A_g}|T_g f(x)|^{n\over n-\alpha}\bigg]d\nu(x)\le C(1+|F|)\big(1+\log^+|F|+\nu(E)\big).\eqdef{18}$$
For given $E$ and $F$, if $g$ is H\" olderian of order $\delta\in(0,1]$ on $S^{n-1}$ and there exists a ball $B(x_0,r_0)$ such that $|B(x_0,r_0)\cap F|=|B(x_0,r_0)|$, and $\nu(B(x_0,r)\cap E)\ge c_1 r^{\sigma n }$ for $r\le r_0$, with $c_1>0$,  then the 
exponential constant is sharp.
\smallskip
b)  If $\nu=$Lebesgue measure, then \eqref{18} can be improved as 
$$\int_E \exp\bigg[{1\over A_g}|T_g f(x)|^{n\over n-\alpha}\bigg]dx\le C(|E|+|F|\big),\eqdef{18imp}$$
and the inequality is sharp under the same conditions as in a).
\par

The condition that $F$  contains a  ball (up to a set of zero measure) which has enough mass shared by $E$, is essentially necessary in order to guarantee sharpness in the above corollary. In general the sharp exponential  constant will depend on the relative geometry of the sets $E$ and $F$: the less  the mass they have in common, the larger  the sharp constant. This is a reflection of the fact that the potential becomes ``less effective" as the sets $E$ and $F$ get  more and more separated (in this regard, see [FM1, Remark 3, p. 5112]). 

Note that \eqref{18imp} in the case $E=F$ vas first obtained by Cohn-Lu [CL, Theorem 1.5].
\bigskip
Finally, we believe that the techniques developed in this paper could be adapted to other settings, such as the Heisenberg group or other noncompact manifolds. For example, on the Heisenberg group one could consider a version of Theorem 1 for the powers of the sublaplacian, which would extend the results in [LL1] to higher order operators.

\vskip1em

\centerline{\scaps 2. Adams inequalities for general homogeneous Riesz-like potentials on $\R^n$}\bigskip
From now on $g:\Rn\setminus\{0\}\to\R$ will denote a homogeneous function of order $\a-n$, that is 
$$g(x)=g(x^*) |x|^{\a-n},\qquad x^*={x\over|x|}\in S^{n-1}$$
and $T_g$ will denote the convolution operator 
$$T_g f(x)=(g* f)(x)=\irn g(x-y)f(y)dy.\eqdef{Tg}$$
We define  $T_g f$ on vector-valued functions in the same way, with the understanding that 
$$g=(g_1,...,g_m),\quad f=(f_1,...,f_m),\quad gf=g_1f_1+...+g_mf_m,\quad |f|=(f_1^2+...+f_m^2)^{1/2},\eqdef{11b}$$
where each $g_j$ is homogeneous of order $\a-n$.
Additionally, we will let, for $f=(f_1,..,f_m)$
$$f^+=(f_1^+,....,f_m^+),\qquad f^-=(f_1^-,...,f_m^-)\eqdef{11c}$$
where $f_j^+$ and $f_j^-$ denote the positive and negative parts of $f_j$. We will say that a vector $f$ is nonnegative if each component of $f$ is nonnegative. The results and their proofs below are valid for the vector-valued case with the above conventions;  we will not distinguish between the scalar case and the vector-valued case, except in a few isolated instances.

\proclaim Theorem 7 (Adams inequality \`a la Ruf).
If $0<\alpha<n$ and $g$ is nonnegative and Lipischitz, then there  
there exists a constant $C=C(\alpha,g,n)$ such that for every measurable $E\subseteq \R^n$ with $|E|<\infty $ and for all compactly supported $f$ with
$$\|f\|_{\na}^\na+\|T_g f\|_{n/\a}^{n/\a}\le 1\eqdef {11}$$
we have
$$\int_E\exp\bigg[{1\over A_g}\,|T_g f(x)|^{n\over n-\a}\bigg]dx\le C(1+|E|).\eqdef {12}$$
where $$A_g={1\over n}\isn |g(\omega)|^\nna d\omega,\eqdef A$$
and also
$$\int_{\sR^n}\exp_{[\nsa-2]}\bigg[{1\over A_g}\,|T_g f(x)|^{n\over n-\a}\bigg]dx\le C.\eqdef {12aa}$$
If $g\in C^{2n}(S^n)$ and  $|E|>0$, then the exponential constant $A_g^{-1}$ in  \eqref {12}, and also in \eqref{12aa}, is sharp, i.e. it cannot be replaced by a larger number.\hfill\break \indent
 If $g$ changes sign then the above results continue to hold for all compactly supported $f$ satisfying the additional pointwise condition: for all $a\in \R^n$
$$|T_g f(x)|\le \int_{|y-a|\le 2}|g(x-y)|\,|f(y)|dy+C_1\|T_g f\|_\na,\qquad 
\;|x-a|\le 1\eqdef{12q}$$
almost everywhere, where $C_1$ is a constant depending only on $\alpha$ and $n$.
\par 
\bigskip\def\s{\sigma}
It is important to point out that the constant $C$ on the right-hand side in  \eqref{12} is independent of the measure  of $\supp f$. Indeed, as stated in Corollary 6,  there exists $C$ depending only on $\alpha,n,g$ so that  for  given measurable sets $E,F$ with finite measure and for all  functions $f$ with $\supp f\subseteq F$ and  $\| f\|_\na\le1$ we have 
$$\int_E\exp\bigg[{1\over A_g}\,|T_g f(x)|^{n\over n-\a}\bigg]dx\le C\big(|E|+|F|\big)\eqdef {12h}$$
(the nonnegativity of $g$ is not needed here). We will give a proof later, however this  inequality can first be verified for $|F|=1$ using O'Neil's inequality and  the Adams-Garsia lemma as in [FM1], followed by a dilation.
The point is that one cannot hope to make the right hand side of  \eqref{12h} to be independent of $|F|$ without further restrictions on $\|T_g f\|_\na$. Indeed, it is possible to find a family of  functions $f_r$ supported on balls $F=B(0,r)$ with   $\|f_r\|_\na\le~1$, and with even the further restriction that $\|T_g f_r\|_\na\le 1$,  for which the inequality in  \eqref{12h} is reversed for any given $E$ with positive measure, for all $r$ large enough (see proof of Corollary 8).
The stronger ``Ruf condition" in \eqref{11} is precisely  what it is needed in order to compensate for the lack of control on the support of $f$, or rather its measure.

In the same spirit as Corollary 2 we have the following general Adachi-Tanaka type result:

\proclaim Corollary 8 (Adams inequalities \`a la Adachi-Tanaka). If $0<\alpha<n$  and $g$ is nonnegative and Lipschitz, then there exists $C=C(\a,g,n)$ such that for $0<\theta<1$, $1<q\le+\infty$,   for all   $E\subseteq  \Rn$ with $|E|<\infty$ 
and for all compactly supported $f$ with
$$\cases{\|f\|_\na^{q\na}+\| T_g f \|_\na^{q\na}\le 1, &  if $q<\infty$\cr\cr   \max\big\{\| f\|_\na,\,\|T_g f\|_\na\big\}\le 1, & if  $q=+\infty.$\cr}\eqdef{103}$$ we have
$$\int_E \exp\Big[{\theta\over A_g}|T_g f(x)|^\nna\Big]dx\le C( 1-\theta)^{-{1\over q'}}(1+|E|)\eqdef{102}$$
and also 
$$\int_{\sR^n}\exp_{[\nsa-2]}\Big[{\theta\over A_g}|T_g f(x)|^\nna\Big]dx\le C( 1-\theta)^{-{1\over q'}}.\eqdef{102aa}$$
Inequalities  \eqref{102}  and \eqref{102aa} are sharp, in the sense that if $|E|>0$ and $g\in C^{2n}(S^n)$ then   for given $q\in(1,\infty]$,  there exists a family of functions $f_\theta\in L^\nsa(\Rn)$ satisfying \eqref{103} for which the inequalities \eqref{102} and \eqref{102aa} are reversed; in particular, the exponential integral cannot be uniformly bounded if $\theta=1$.  
\par

We observe here that the proof of the original Adachi-Tanaka estimate ($q=+\infty$) without sharp control on the right-hand side in terms of $\theta$, does not require the full Theorem~1, and it is much simpler to prove directly with the methods given in this paper. We will present this argument later in this section.

\bigskip

\centerline{\bf Proof of Theorem 7: Overview and preliminary lemmas}\bigskip

Let us first state an elementary lemma, which is more like an observation, with the hope of  clarifying once and for all the equivalence between exponential inequalities on sets of finite measure and regularized exponential inequalities on sets of arbitrary measure.

\smallskip\proclaim Lemma 9 (Exponential Regularization Lemma). Let  $(N,\nu)$ be a measure space and $1<p<\infty$, $\a>0$. Then for every $u\in L^p(N)$ we have
$$\int_{\{u\ge1\}} e^{\alpha|u|^{p'}}d\nu-e^\alpha\|u\|_p^p\le\int_N\bigg(e^{\a |u|^{p'}}-\sum_{k=0}^{[p-2]}{\a^k |u|^{kp'}\over k!}\bigg)d\nu\le\int_{\{u\ge1\}} e^{\alpha|u|^{p'}}d\nu+e^\a\|u\|_p^p.$$
In particular, the functional $\int_N \exp_{[p-2]}^{}\big[\a |u|^{p'}\big]$ is  bounded on a bounded subset $X$ of  $L^p$, if and only if $\int_{\{u\ge 1\}}\exp\big[\a |u|^{p'}\big]$ is  bounded  on $X$.
\par
\pf Proof. Recall that $[p-2]$ is the smallest integer greater or equal $p-2$ . To start, write $N=(N\cap\{u<1\})\cup (N\cap\{u\ge1\})$, and split the middle integral accordingly.  The estimates then follow from  the Taylor's series of $e^{\a |u|^{p'}}$ and straightforward considerations.\endpf
Note also that for any measurable $E\subseteq N$ with $\nu(E)<\infty$  we obviously have
$$\int_E e^{\alpha|u|^{p'}}d\nu\le \int_{\{u\ge1\}} e^{\alpha|u|^{p'}}d\nu+e^\a\nu(E).$$

 From now on we will only focus on   exponential integrals over sets of finite measure,
since all the inequalities that involve regularized exponentials  stated in this paper can be deduced at once from this case, just by appealing to the lemma above.

 Let us now note that it is enough to prove the theorem in  the case $|E|\le 1$. This is because of a dilation argument that will be also used later on. 
The observation  is that if  $f_\lambda(x)=\lambda^\a f(\lambda x)$ then 
$$\|f_\lambda\|_\na=\|f\|_\na,\qquad T_gf_\lambda(x/\lambda)=T_g f(x),\qquad \|T_g f_\lambda\|_\na^\na=\lambda^{-n}\|T_g f\|_\na^\na.\eqdef{D1}$$
and 
$$\int_E \exp\Big[{1\over A_g}|T_g f(x)|^\nna\Big]dx=\lambda^n\int_{E/\lambda} \exp\Big[{1\over A_g}|T_g (f_\lambda)(x)|^\nna\Big]dx.\eqdef{D2}$$
Hence, \eqref{12h} follows for general $E$ with finite measure, if it's known for $|E|=1$, by taking $\lambda=|E|^{1/n}$.

To get to the heart of the matter, we now state and prove a  Lemma, which is very elementary in nature, yet crucial  in the proof of Theorem 7 and other results of this paper:\medskip

\proclaim Lemma 10. Let  $(N,\nu)$ be a measure space, and let $V,Z$ be  vector spaces of measurable functions (real, complex or vector valued) on a measurable set $E\subseteq N$. Let $T:V\to Z$ be an operator such that $T(\lambda f)=\lambda Tf$ for any $f\in V$ and any $\lambda\ge~0$,  and let $p:V\to [0,\infty]$ be a seminorm. Finally let $\b>1$ and $\b'=\ds{\b\over\b-1}$. \hfill\break 
\indent If there exists $c_0$ such that for a fixed subset $V_0\subseteq V$
$$\int_E \exp\bigg[{1\over A}|T f(x)|^{\beta}\bigg]d\nu(x)\le c_0,\qquad \forall f\in V_0,\, p(f)\le 1$$
then, for each $K>0$\def\bbp{{\b\over\bp}} and for each $f\in V_0$ with $p(f)\le 1$ we have 
$$(A)\qquad\quad  \int_E \exp\bigg[{1\over A}\big(|T f(x)|+K\big)^{\beta}\bigg]d\nu(x)\le c_0 \exp\bigg[{1\over A} \bigg({K^{\bp}\over1-p(f)^\bp}\bigg)^{\bbp}\,\bigg]$$
and 
$$(B)\qquad \qquad \int_E \exp\bigg[{1\over A}\Big(|T f(x)|+K\big(1-p(f)^\bp\big)^{1\over\bp}\Big)^\b\,\bigg]d\nu(x)\le c_0e^{ {1\over A}K^\b}.$$
\par

  \eject 

\pf Proof. From H\"older's inequality we have 
$$a\theta^{\bpi}+b(1-\theta)^\bpi\le (a^\b+b^\b)^\bi,\qquad a,b\ge0,\,0\le\theta\le 1,\eqdef{20x}$$
from which estimate (B) follows, with $\theta=p(f)^\bp>0$, $a=T\big(f/ p(f)\big),\, b=K$. Clearly (A) is just another way of writing (B), since $K$ can be arbitrary.\endpf 
 
   \smallskip

As it is apparent from the proof, there is nothing peculiar about exponential integrability in this lemma, and (A) and (B) are clearly equivalent. Nonetheless, we find it convenient to have the estimates  in  (A) and (B) explicitly stated as above, since they will be used directly several times.
The first main application of the Lemma is in Theorem~7,  with $V=\{f\in L^\nsa(\Rn),\,\supp f $ compact$\}$, $Z=\{u:\R^n\to\R^m\,\,{\hbox{a.e. finite}}\}$, $\b=n/\a$, $p(f)=\|f\|_\na$, and $E\subseteq\Rn$ with $\nu=$Lebesgue measure (and $|E|\le 1$), and obviously $T=T_g$.

For the benefit of the reader we will now summarize the main strategy behind the proof of Theorem 1. We start from the Adams inequality
$$\int_E\exp\bigg[{1\over A_g}\,|T_g f_1(x)|^{n\over n-\a}\bigg]dx\le C$$
valid for all functions $f_1$ with $\|f_1\|_\na\le1$ with $|\supp f_1|\le \kappa$, with $\kappa$ depending  possibly  only on $n$ (this follows from \eqref{12h}). The main idea is that the inequality in Theorem 7 is true if we can  write $T_g f=T_g f_1+T_g f_2$, where $f=f_1+f_2$, $\,|\supp f_1|\le \kappa$, and where the additive perturbation $T_g f_2$ satisfies either
$$a)\;\;|T_g f_2(x)|\le C, \quad {\hbox {if}}\quad  x\in  E\qquad {\hbox{and}}\quad \|f_1\|_\na\le \theta_n<1$$
or
$$b)\;\;|T_g f_2(x)|\le C\big(1-\|f_1\|_\na^\na\big)^\asn, \quad {\hbox{if}}\quad  x\in E\qquad {\hbox{and}}\quad 0<\theta_n\le \|f_1\|_\na\le1, $$
 $\theta_n$ being a suitable explicit constant depending only on $n$. It is clear that in either case a) or b) one can apply (A) or (B) respectively, to derive the desired inequality.   The original set $E$ and the original function $f$  will  be suitably split so as to reduce matters to estimates a) and b). To this end, we will consider several scenarios  depending on where the $L^\na$ masses of $f,\, f^+,\,f^-$  are concentrated, and on where the potential is pointwise positive. An estimate for $T_g f_2$ as in a) will follow if the  function $f_2$ is either pointwise  small, or if its  support is ``well separated" from $E$; this will be a consequence of Lipschitz estimates for $T_g f$. An estimate as in b) will instead occur, roughly speaking,  when both $E$ and the mass of $f$ are concentrated in a fixed ball, and it will be the most critical case of the proof, the only one where Ruf's condition is needed.

\eject

We now establish some regularity estimates for the operator $T_g$. 

\proclaim Lemma 11. Let $ f\in L^{\nsa}(\Rn)$, compactly supported and $F\subseteq \Rn$ a closed set such that  either
\smallskip\noin
(i) $|\supp f|\le 1$ and dist$ (F,  \supp f)\ge R\ge1$
\smallskip
\noin or
\smallskip\noin
(ii) $\,\supp f\subseteq B_{2R}^c$, $F\subseteq B_{R}$, 
\smallskip\noin
then $T_g f$ is Lipschitz on $F$, in particular there exists $D=D(n,\a)$ such that 
$$|T_g f(x_1)-T_g f(x_2)|\le{D\over R}\| f\|_\na |x_1-x_2|,\qquad x_1,x_2\in F.\eqdef{21}$$
Moreover, if $x_0\in F$ is so that $\max_{x\in F} |T_g  f(x)|=|T_g  f(x_0)|,$ then for $R\ge 2$
$$\max_{x\in F} |T_g  f(x)|\le \bigg(\avg_{B_1(x_0)}|T_g  f(x)|^\nsa dx\bigg)^{\asn}+{2D\over R}\| f\|_\na.\eqdef{22}$$\par
\pf Proof. We have, for all $y,x_1,x_2\in\R^n,\, y\neq x_1,\,y\neq x_2$ 
$$\big||x_1-y|^{\a-n}-|x_2-y|^{\a-n}\big|\le|x_1-x_2|(n-\a)\Big(\min\big\{|x_1-y|,\,|x_2-y|\big\}\Big)^{\a-n-1}\eqdef{22a}$$
and since $g$ is Lipschitz on the sphere (in terms of the Euclidean distance) we find that  there exists $L=L(g,\a,n)$ such that
$$|g(x_1-y)-g(x_2-y)|\le L|x_1-x_2|\Big(\min\big\{|x_1-y|,\,|x_2-y|\big\}\Big)^{\a-n-1}.\eqdef{22ab}$$
Hence, for $x_1,x_2\in F$ and each $y\in \supp f$
$$\big|g(x_1-y)-g(x_2-y)\big|\le L |x_1-x_2|\cases {R^{\a-n-1} & in case (i)\cr
2^{n+1-\a}|y|^{\a-n-1} & in case (ii).\cr}$$
and  in case (i) 
$$|T_g f(x_1)-T_g f(x_2)|\le L|x_1-x_2|R^{\a-n-1}\int_{\supp f}|f(y)|dy\le {D\over R}\|f\|_\na |x_1-x_2|,$$
whereas  in case (ii)
$$\eqalignno{&|T_g f(x_1)-T_g f(x_2)|\le L|x_1-x_2|2^{n+1-\a}\int_{B_{2R}^c}|f(y)||y|^{\a-n-1}dy & \eqdef{22k}\cr&\le L |x_1-x_2|
2^{n+1-\a}\|f\|_{\na}\bigg(\int_{B_{2R}^c}|y|^{-n-\nna}dy\bigg)^{{n-\a\over n}}\le{D\over R}\|f\|_\na |x_1-x_2|.\cr}$$
Note that since $\supp f$ is compact we have $|T_g f(x)|\le C_f |x|^{\a-n}$ for large $|x|$ (where $C_f$ depends on $f,\a,n, g$), hence the supremum of $|T_g f|$ on $F$ is attained in $F$ in either  case (i) or (ii). Let $M=\max_{x\in F} |T_g f(x)|=|T_g f(x_0)|$, some $x_0\in F$. After a moment's reflection the reader should realize that if $R\ge2$ then estimate \eqref{21} holds also for all $x_1,x_2$ in the set $F\cup B_1(x_0)$ (whose distance from $\supp f$ is at most $R-1$), by possibly enlarging slightly the constant~$D$.
 If $M\le {2D\over R}\|f\|_\na$ then \eqref{22} is true. If $M>{2D\over R}\|f\|_\na$ then  using \eqref{21} for  $x_1,x_2\in F\cup B_1(x_0)$ we get $|T_g f(x)|\ge M-{D\over R}\|f\|_\na\ge0$ for $x\in B_1(x_0)$ and 
$$\bigg( M-{D\over R}\|f\|_\na\bigg)^\nsa |B_1|\le \int_{B_1(x_0)} |T_g f(x)|^\nsa dx$$
which yields \eqref{22}.\endpf

\bigskip \proclaim Lemma 12. If $f$ is compactly supported and $|f|\le 1$ on $\Rn$, then  there exists $D=D(n,\a)$ such that
$$|T_g f(x_1)-T_g f(x_2)|\le D(1+\|f\|_\na)\min\big\{|x_1-x_2|^\a,|x_1-x_2|\big\},\,\qquad x_1,x_2\in \Rn.\eqdef{23}$$
If $\alpha=1$ there exists $D=D(n)$ so that
$$|T_1 f(x_1)-T_1 f(x_2)|\le D(1+\|f\|_n)|x_1-x_2|\Big(1+\log^+{1\over|x_1-x_2|}\Big),\,\qquad x_1,x_2\in \Rn,\;x_1\neq x_2.\eqdef{24}$$
\smallskip If $x_0\in\Rn$ is so that $\max_{x\in \sR^n} |T_g f(x)|=|T_g f(x_0)|$ then
$$\|T_g f\|_\infty\le \bigg(\avg_{B_1(x_0)}|T_g f(x)|^\nsa dx\bigg)^{\asn}+{2D}(1+\|f\|_\na).\eqdef{25}$$\par

 \par
\pf Proof. For each $x_1,x_2\in\Rn$ we have
$$\eqalign{&T_g f(x_1)-T_g f(x_2)= \!\!\!\!\mathop\int\limits_{|y-x_1|\le {1\over3}|x_1-x_2|}\hskip-2em g(x_1-y)f(y)dy+\hskip-1em \mathop\int\limits_{|y-x_2|\le {1\over3}|x_1-x_2|}\hskip-2emg(x_2-y)f(y)dy\cr&+ \hskip-1em\mathop\int\limits_{{|x_1-y|>|x_2-y|\atop |y-x_2|> {1\over3}|x_1-x_2|}}\hskip-2em\Big(g(x_2-y)-g(x_1-y)\Big)f(y)dy+\hskip-1em\mathop\int\limits_{{|x_2-y|>|x_1-y|\atop |y-x_1|> {1\over3}|x_1-x_2|}}\hskip-2em\Big(g(x_1-y)-g(x_2-y)\Big)f(y)dy.\cr}$$
then,  using \eqref{22ab} and $|g(x)|\le C|x|^{\a-n}$ we get
$$\eqalign{|T_g f(x_1)&-T_g f(x_2)|\le 2C\omega_{n-1}\int_0^{ {1\over3}|x_1-x_2|}r^{\a-1}dr+2L|x_1-x_2|\hskip-3em\mathop\int\limits_{{|y-x_1|<y-x_2|\atop  {1\over3}|x_1-x_2|<|y-x_1|\le 1}}\hskip-2em |x_1-y|^{\a-n-1}dy\cr&\hskip12em+2L|x_1-x_2|\mathop\int\limits_{{|y-x_1|\ge 1}} |x_1-y|^{\a-n-1}|f(y)|dy\cr&\le
{2C\omega_{n-1}\over \a 3^\a}|x_1-x_2|^\a+2\omega_{n-1}L|x_1-x_2|\int_{{1\over3}|x_1-x_2|}^1 r^{\a-2}dr\cr&\hskip12em +2L\Big({n-\a\over n}\Big)^{\asn}\|f\|_\na|x_1-x_2|\cr
}$$
and this proves \eqref{23} and \eqref{24}. Clearly $T_g f$ is continuous on $\Rn$ and $|T_g f(x)|\le C|x|^{a-n}$ for large $|x|$, so  $|T_g f(x)|$ has a maximum at some $x_0$. Estimate \eqref{25} is obtained as in Lemma 11.\endpf
\vskip2em
\proclaim Lemma 13. If $f\in L^{\nsa}(\Rn)$, $\,|\supp f|\le \kappa$, then there is $C=C(n,\a,g)$ such that  for all $E\subseteq \Rn$ with $|E|<\infty$
$$\bigg(\int_E |T_g f|^{\nsa}\bigg)^\asn\le C(|E|^{\asn}+\kappa^\asn) \|f\|_\na.\eqdef{25a}$$
\par\pf Proof. This  follows at once from O'Neil's inequality: for $t\le|E|$
$$(T_g f)^{**}(t)\le \nsa A_g^{n-\a\over n} t^\asn f^{**}(t)+A_g^{n-\a\over n} \int_t^\kappa s^{\asn-1}f^*(s)ds\le C(|E|^\asn+\kappa^\asn) f^{**}(t).$$
\endpf
\bigskip
Let now $f$ be compactly supported and such that 
$$\|f\|_\na^\na\le 1,\qquad \|T_g f\|_\na^\na\le1.\eqdef{25aa}$$
From now on we let

$$f=f_\ell+f_s,\qquad \quad f_\ell(x)=\cases{f(x) & if $|f(x)|\ge 1$\cr 0 & if $|f(x)|<1$.\cr}\eqdef{26}$$
Obviously   $\supp f_\ell$ is compact, $|f_\ell|\ge1$ in $\supp f_\ell$ and $|\supp f_\ell|\le1.$ Also, $|f_s|\le 1$ on $\Rn$, and $\supp f_s$ is compact, with no control on its measure.
\bigskip
\proclaim Lemma 14. Under condition \eqref{25aa} we have
$$\|T_g f_s\|_{L^\infty(\sR^n)}\le C.\eqdef{27}$$
If $|F|<\infty$ then 
$$\|T_g\big( f_s\chi_F\big)\|_{L^\infty(\sR^n)}\le C(1+|F|^{\asn}).\eqdef{27a}$$
and if $|F^c|<\infty$ then 
$$\|T_g\big( f_s\chi_F\big)\|_{L^\infty(\sR^n)}\le C(1+|F^c|^{\asn}).\eqdef{27b}$$
\par
\pf Proof. From Lemma 12 we have
$$\|T_g f_s\|_\infty \le 2D(1+\|f_s\|_\na)+\bigg(\avg_{B_1(x_0)} |T_g f_s|^\nsa\bigg)^{\asn}$$
where $x_0$ is a  maximum for $|T_g f_s|$.
By Lemma 13
$$\bigg(\int_{B_1(x_0)} |T_g f_s|^\nsa\bigg)^{\asn}\le \bigg(\int_{B_1(x_0)} |T_g f|^\nsa\bigg)^{\asn}+\bigg(\int_{B_1(x_0)} |T_g f_\ell|^\nsa\bigg)^{\asn}\le 1+C\|f_\ell\|_\na\le C.$$
If $|F|<\infty$ then \eqref{27a} follows from \eqref{25}, \eqref{25a}, \eqref{25aa}. If $|F^c|<\infty$ then
just write $f_s\chi_F=f_s-f_s\chi_{F^c}$.
\endpf

\bigskip \centerline{\bf  Proof of Theorem 7.}\bigskip

Let $E\subseteq \R^n$ be measurable and with $|E|\le 1,$ and let $f$ be compactly supported and satisfying \eqref{25aa}. 

After  suitable translations we can  assume that
$$\int_{x_i\le0} |f(x)|^\nsa dx=\int_{x_i\ge0}|f(x)|^\nsa dx={1\over2}\; \nf.\eqdef{29}$$
Define the $2n$ half-spaces
$$H_i^+=\{x\in\Rn: \;x_i\ge4\},\qquad H_i^-=\{x\in\Rn: \; x_i \le -4\}\eqdef{30}$$
\medskip
We organize the proof in 6  cases:
\medskip
\noin{\bf \underbar{Case 1}:} There exists a half-space $H\in \{H_i^+,H_i^-\}$ such that 
$$\int_H |f(x)|^\nsa dx\ge {1\over4n}\nf.\eqdef{31}$$
{\it (The mass of $f$ is not concentrated.)}
\medskip
\noin{\bf \underbar{Case 2}:} 
$$\int_{B_{4\sqrt n}}|f(x)|^\nsa dx\ge{1\over2}\nf,\quad{\hbox {and}}\quad \|f_\ell\|_\na^\na\le{3\over4}\nf.\eqdef{32}$$({\it The mass of $f$ is concentrated but the one of $f_\ell$  is too small.})
\medskip
\noin{\bf \underbar{Case $\bf {2^+}$}:}  The kernel $g$ is {\it {nonnegative}} and 
$$\int_{B_{4\sqrt n}}|f(x)|^\nsa dx\ge{1\over2}\nf,\quad \|f_\ell^+\|_\na^\na\le{3\over4}\nf,\;\;\|f_\ell^-\|_\na^\na\le{3\over4}\nf\eqdef{32a}$$
({\it The mass of $f$ is concentrated but the ones of $f_\ell^+$ and $f_\ell^-$ are too small.})

\medskip
\noin{\bf \underbar{Case 3}:}     $E\subseteq B_{16\sqrt n}^c$ and 
$$\int_{B_{4\sqrt n}}|f(x)|^\nsa dx\ge{1\over2}\nf,\quad \|f_\ell\|_\na^\na\ge{3\over4}\nf.\eqdef{32aa}$$

\noin({\it The masses of $f$ and $f_\ell$ are concentrated but $E$ is too far away.})

\medskip\noin{\bf \underbar{Case $\bf 4$}:}  $E\subseteq B_{16\sqrt n}$ and conditions \eqref{11}, \eqref{12q} hold.\medskip

\noin ({\it Ruf's condition combined with the strong pointwise condition in \eqref{12q} do everything in this case}.) 

\medskip\noin{\bf \underbar{Case $\bf 4^+$}:}   $E\subseteq B_{16\sqrt n}$, the kernel $g$ is {\it {nonnegative}} and 
$$\int_{B_{4\sqrt n}}|f(x)|^\nsa dx\ge{1\over2}\nf,\quad \|f_\ell^+\|_\na^\na\ge{3\over4}\nf\;\;{\hbox{ OR}}  \;\;\|f_\ell^-\|_\na^\na\ge{3\over4}\nf.\eqdef{32b}$$
({\it The mass of $f$ and the mass of $f_\ell^+$, or $f_\ell^-$, are both concentrated near $E$}).\smallskip

Within this case we will consider the  following subcases:\medskip
\smallskip\item{\bf{(i)}} $T_g f(x)\le0$ for any $x\in E$;\smallskip
\item{\bf{(ii)}}  $T_g(f\chi_{B_{32\sqrt n}}^{})(x)$ and $T_g(f\chi_{B_{32\sqrt n}^c}^{})(x)$ have opposite sign, and $T_gf(x)\ge0$, for any $x\in E$;\smallskip
\item{\bf{(iii)}}  $T_g(f\chi_{B_{32\sqrt n}}^{})(x)\ge0$ and $T_g(f\chi_{B_{32\sqrt n}^c}^{})(x)\ge0$ for any $x\in E$, and condition \eqref{11} holds.\smallskip

\medskip
It is easy to see that once  Theorem 7 is proved in the above cases then it is proved in full generality. Indeed,  if Case 1 is not verified then $\int_{B_4\sqrt n}|f|^{\nsa}\ge \half\nf$, so that from Cases 1,2,3 the theorem follows when $E\subseteq B_{16\sqrt n}^c$. For arbitrary $E$ write $E=(E\cap  B_{16\sqrt n})\cup(E\cap  B_{16\sqrt n}^c)$ and the theorem follows if \eqref{12q} is assumed. Similarly the theorem follows from 1, $2^+$, $3^+$, $4^+$ if $g$ is nonnegative.

It is worth emphasizing  that  {\it cases 1, 2, 3 and 4  hold for vector-valued kernels with arbitrary sign}, and that {\it the stronger ``Ruf condition" \eqref{11} is only used  in cases $4$ and $4^+$ }(iii).

We will now prove the main estimate \eqref{12}  in these cases.
\medskip
\noin\underbar {\sl Proof of \eqref{12} in Case 1.}  Suppose  WLOG that 
$$\int_{x_1\ge4} |f(x)|^\nsa dx\ge{1\over4n}\nf\eqdef{33}$$
and write
$$\int_E\exp\bigg[{1\over A_g}\,|T_g f(x)|^{n\over n-\a}\bigg]dx=\int_{E\cap\{x_1\le 2\}}+\int_{E\cap\{x_1\ge2\}}=I+II.\eqdef {34}$$
To estimate $I$ write
$$T_g f=T_g\big(f_\ell\chi_{\{x_1\le 4\}}^{}\big)+T_g\big(f_\ell\chi_{\{x_1\ge 4\}}^{}\big)+T_g f_s.\eqdef{35}$$
From Lemmas 11 and 13 we get 
$$\big|T_g\big(f_\ell\chi_{\{x_1\ge 4\}}^{}\big)(x)\big|\le \bigg(\avg_{B_1(x_0)}\big|T_g\big(f_\ell\chi_{\{x_1\ge 4\}}^{}\big) (x)\big|^\nsa dx\bigg)^{\asn}+2D\|f_\ell\|_\na\le C,\quad  
\forall x\in\{x_1\le 2\}$$
(here $x_0$ is a maximum for $\big|T_g\big(f_\ell\chi_{\{x_1\ge 4\}}^{}\big) (x)\big|$ in $\{x_1\le 2\}$.)
From Lemma 14 we then have
$$|T_g f(x)|\le \big|T_g\big(f_\ell\chi_{\{x_1\le 4\}}^{}\big) (x)\big|+K,\qquad \forall x\in\{x_1\le 2\}$$
some $K$ depending only on $\a,n$, and (A) of Lemma 10 implies that $I\le C$, since $$\|f_\ell\chi_{\{x_1\le 4\}}^{}\|_\na^\na\le\|f\chi_{\{x_1\le 4\}}^{}\|_\na^\na\le \Big(1-{1\over 4n}\Big)\nf<1-{1\over 4n}<1.$$
The estimate of $II$ is similar, this time write
$$T_g f=T_g\big(f_\ell\chi_{\{x_1\le 0\}}^{}\big)+T_g\big(f_\ell\chi_{\{x_1\ge 0\}}^{}\big)+T_g f_s.\eqdef{36}$$
and Lemmas 11,  13, 14, imply
$$|T_g f(x)|\le \big|T_g\big(f_\ell\chi_{\{x_1\ge 0\}}^{}\big) (x)\big|+K,\qquad \forall x\in\{x_1\ge 2\}$$
some $K$ depending only on $\a,n$, and (A) of Lemma 10 implies that $II\le C$, since 
$$\|f_\ell\chi_{\{x_1\ge0\}}^{}\|_\na^\na\le\|f\chi_{\{x_1\ge 0\}}^{}\|_\na^\na= {1\over2}\nf<{1\over2}.$$
\medskip

\eject

\noin\underbar {\sl Proof of \eqref{12} in Case 2.} \smallskip
Assume \eqref{32} and write
$$\int_E\exp\bigg[{1\over A_g}|T_g f(x)|^\nna\bigg] dx\le\int_{E}\exp\bigg[{1\over A_g}\Big(T_{|g|} |f_\ell|(x)+|T_g f_s(x)|\Big)^\nna \bigg]dx$$
 and \eqref{12} follows from \eqref{27a} and (A) of Lemma 10 applied to $f_\ell$, since $\|f_\ell\|_\na^\na\le {3\over4}$.
\bigskip
\noin\underbar {\sl Proof of \eqref{12} in Case $2^+$.} \smallskip
Assume \eqref{32a} and  write
$$\eqalign{\int_E\exp\bigg[{1\over A_g}|&T_g f(x)|^\nna\bigg] dx\le\int_{E^+}\exp\bigg[{1\over A_g}\Big(T_g f_\ell^+(x)+|T_g f_s(x)|\Big)^\nna \bigg]dx+\cr&+\int_{E^-}\exp\bigg[{1\over A_g}\Big(T_g f_\ell^-(x)+|T_g f_s(x)|\Big)^\nna\bigg] dx\cr}$$
where $E^{\pm}=\{x\in E:\,T_g f_\ell^\pm(x)\ge T_g f_\ell^\mp(x)\}$, and \eqref{12} follows from (A), applied to $f_\ell^\pm$.

\bigskip
\noin\underbar {\sl Proof of \eqref{12} in Case 3.} Suppose $E\subseteq  B_{16\sqrt n}^c$ and that the estimates in  \eqref{32aa} hold.  
For $x\in E$ we then have 
$$|T_g f(x)|\le \big|T_g\big(f_\ell\chi_{B_{8\sqrt n}}^{}\big)(x)\big|+\big
|T_g\big(f_\ell\chi_{B_{8\sqrt n}^c}^{}\big)(x)\big|+|T_g f_s(x)|\le\big|T_g\big(f_\ell\chi_{B_{8\sqrt n}^c}^{}\big)(x)\big|+ C$$
from Lemmas 11-14.
  On the other hand, since $\big\|f\chi_{B_{4\sqrt n}}^{}\big\|_\na^\na\ge{1\over2}\nf$, it must be that $\big\|f_\ell\chi_{B_{8\sqrt n}^c}^{}\big\|_\na^\na\le{1\over2}\nf<{1\over2}$, and $\big|\supp\big(f_\ell\chi_{B_{8\sqrt n}^c}^{}\big)\big|\le 1$, so that our estimate \eqref{12}  follows again from (A) of Lemma 10,  applied  to $f_\ell\chi_{B_{8\sqrt n}^c}^{}$.
\bigskip
\noin\underbar {\sl Proof of \eqref{12} in Case $4$.} \smallskip
Let  $E\subseteq  B_{16\sqrt n}$ and suppose that $f$ satisfies both  \eqref{11}  and \eqref{12q}. Estimate \eqref{12q} implies that 
$$|T_g f(x)|\le \int_{B_{32\sqrt n}}|g(x-y)|\,|f(y)|dy+C_1\|T_g f\|_\na,\qquad  x\in B_{16\sqrt n}.$$
(if $x\in B_{16\sqrt n}$, pick $a\in B_{16\sqrt n}$ such that $|x-a|\le 1$, use \eqref{12q}, and enlarge the domain of integration.)

The Ruf condition $\nf+\ntf\le1$  gives
$$ \eqalign{|T_gf(x)|&\le T_{|g|}^{}\big(|f|\chi_{B_{32\sqrt n}}^{}\big)+C(1-\|f\|_\na^\na)^\an\cr& \le  T_{|g|}^{}\big(|f|\chi_{B_{32\sqrt n}}^{}\big)+C(1-\|f\chi_{B_{32\sqrt n}}^{}\|_\na^\na)^\an\cr}\qquad x\in B_{16\sqrt n,}\eqdef{usaruf1}$$
and  Lemma 10(B) yields inequality \eqref{12}.
\bigskip
\noin\underbar {\sl Proof of \eqref{12} in Case $4^+$.} \smallskip
It is clear that   in   \eqref{32b} it is  enough to assume that  $\|f_\ell^+\|_\na^\na\ge{3\over4}\nf$, in which case have
$$\|f_\ell^-\|_\na^\na\le {1\over4}\nf,\eqdef{39}$$
and
$$\int_{B_{4\sqrt{n}}} (f_\ell^+)^{\nsa} \ge\irn|f_\ell^+|^\nsa-\int_{B_{4\sqrt n}^c}|f|^{\nsa}\ge{1\over4}\nf.\eqdef{40}$$

In case (i) we have $T_gf(x)\le 0$ in $E$, hence, since $g\ge0$, 
$$\eqalign{\int_E\exp\bigg[{1\over A_g}|&T_g f(x)|^\nna\bigg] dx\le\int_{E}\exp\bigg[{1\over A_g}(T_g f_\ell^-(x)+|T_g f_s(x)|)^\nna\bigg] dx\cr}\le C,$$
using again (A) applied to $f_\ell^-$. 

In case (ii) write
$$\int_E\exp\bigg[{1\over A_g}|T_g f(x)|^\nna\bigg] dx=\int_{E}\exp\bigg[{1\over A_g}
\Big(T_g\big(f\chi_{B_{32\sqrt n}}^{}\big)+T_g\big(f\chi_{B_{32\sqrt n}^c}^{}\big)\Big)^\nna\bigg]dx.$$

If $T_g\big(f\chi_{B_{32\sqrt n}}^{}\big)\ge0 $ and $T_g\big(f\chi_{B_{32\sqrt n}^c}^{}\big)\le0$ then 
$$\int_E\exp\bigg[{1\over A_g}|T_g f(x)|^\nna\bigg] dx\le\int_E\exp\bigg[{1\over A_g}\Big(T_g\big(f\chi_{B_{32\sqrt n}}^{}\big)(x)\Big)^\nna \bigg]dx\le C$$
by the original Adams inequality \eqref{12h}.

If instead $T_g\big(f\chi_{B_{32\sqrt n}}^{}\big)\le0 $ and $T_g\big(f\chi_{B_{32\sqrt n}^c}^{}\big)\ge0$ then 
$$\int_E\exp\bigg[{1\over A_g}|T_g f(x)|^\nna\bigg] dx\le\int_E\exp\bigg[{1\over A_g}\Big(T_g\big(f\chi_{B_{32\sqrt n}^c}^{}\big)(x)\Big)^\nna \bigg]dx\le C$$
since on $E\subseteq B_{16\sqrt n}^{}$ we have
$$0\le T_g\big(f\chi_{B_{32\sqrt n}^c}^{}\big)(x)\le |T_g\big(f_\ell\chi_{B_{32\sqrt n}^c}^{}\big)(x)|+
|T_g\big(f_s\chi_{B_{32\sqrt n}^c}^{}\big)(x)|\le C$$
from lemmas 11-14.

\eject
 
In case (iii), the most critical situation, let us assume the Ruf condition
$\nf+\ntf\le1$
 and write 
$$\eqalign{\int_E\exp\bigg[{1\over A_g}|T_g f(x)|^\nna\bigg] dx\le\int_{E}\exp\bigg[&{1\over A_g}\bigg(T_g\big(f^+\chi_{B_{32\sqrt n}}^{}\big)(x)+\cr&+T_g\Big((f_\ell^++f_s)\chi_{B_{32\sqrt n}^c}^{}\Big)(x)\bigg)^\nna \bigg]dx.\cr}\eqdef{11h}$$

By Lemma 11, for $x\in E\subseteq B_{16\sqrt n}^{}$
$$\eqalign{0\le T_g\Big((f_\ell^++f_s)\chi_{B_{32\sqrt n}^c}^{}\Big)(x)\le\bigg(\avg_{\!\!\!B_1(x_0)} \Big|T_g\Big((f_\ell^++f_s)\chi_{B_{32\sqrt n}^c}^{}&\Big)\Big|^{\nsa}dx\bigg)^{\asn}+\cr&+C\big\|(f_\ell^++f_s)\chi_{B_{32\sqrt n}^c}^{}\big\|_\na\cr}$$
where $x_0$ is a maximum point for $T_g\Big((f_\ell^++f_s)\chi_{B_{32\sqrt n}^c}^{}\Big)$ on 
$\overline{B}_{16\sqrt n}$. 

We have 
$$\eqalign{\bigg(\int_{B_1(x_0)} \Big|T_g\Big((f_\ell^+&+f_s)\chi_{B_{32\sqrt n}^c}^{}\Big)\Big|^{\nsa}dx\bigg)^{\asn}\le\bigg(\int_{B_1(x_0)} \Big(T_g\big(f_\ell^++f_s^++f_s^-\chi_{B_{32\sqrt n}^c}^{}\big)\Big)^{\nsa}dx\bigg)^{\asn}\cr&\le \bigg(\int_{B_1(x_0)} |T_g f|^\nsa dx\bigg)^\asn+\bigg(\int_{B_1(x_0)} \Big|T_g\big(f_\ell^-+f_s^-\chi_{B_{32\sqrt n}}^{}\big)\Big|^{\nsa}dx\bigg)^\asn\cr&\le\|T_g f\|_\na+C\Big(\|f_\ell^-\|_\na+\big\|f_s^-\chi_{B_{32\sqrt n}}^{}\big\|_\na\Big)\cr}$$
(by Lemma 13). Hence, for $x\in E$ 
$$\eqalign{0&\le T_g\Big((f_\ell^++f_s)\chi_{B_{32\sqrt n}^c}^{}\Big)(x)\le\cr& \le C\Big(\|T_g f\|_\na+\|f_\ell^-\|_\na+\big\|f_s^-\chi_{B_{32\sqrt n}}^{}\big\|_\na+\big\|f_\ell^+\chi_{B_{32\sqrt n}^c}^{}\big\|_\na+\big\|f_s\chi_{B_{32\sqrt n}^c}^{}\big\|_\na\Big)\cr&\le C \Big(\|T_g f\|_\na^\na+\|f_\ell^-\|_\na^\na+\big\|f_s^-\chi_{B_{32\sqrt n}}^{}\big\|_\na^\na+\big\|f_\ell^+\chi_{B_{32\sqrt n}^c}^{}\big\|_\na^\na+\big\|f_s\chi_{B_{32\sqrt n}^c}^{}\big\|_\na^\na\Big)^\an\cr&=C\Big(1-\big\|f_\ell^+\chi_{B_{32\sqrt n}}^{}\big\|_\na^\na-\big\|f_s^+\chi_{B_{32\sqrt n}}^{}\big\|_\na^\na\Big)^\an= C\Big(1-\big\|f^+\chi_{B_{32\sqrt n}}^{}\big\|_\na^\na\Big)^\an\cr}\eqdef{usaruf2}$$
where in the second to  last identity we used that $\ntf\le1-\nf$.

By applying (B) of Lemma 10 to the function $f^+\chi_{B_{32\sqrt n}}^{}$ we obtain that the integral in \eqref{11h} is bounded by  constant $C$. This concludes the proof of the inequality part of Theorem~1. 
\bigskip

\eject
  
Now we show that the constant $ A_g^{-1}$ is best possible in \eqref{12}, in the sense that if $|E|>0$ then  we can find a family of functions $\psi_\e\in L^{\na}(\Rn)$ such that $\|\psi_\e\|_\na+\|T_g \psi_\e\|_\na\le 1$ and 
$$\lim_{\e\to0^+}\int_E\exp\bigg[{1+\delta\over A_g}\,|T_g \psi_\epsilon(x)|^{n\over n-\a}\bigg]dx=+\infty,\qquad \forall \delta>0.\eqdef {12x}$$
First notice that for a.e. $x\in E$ we have $|E\cap B_\e(x)|/|B_\e(x)|\to1$, as $\e\to0$, therefore we can assume WLOG that for some $\e_0>0$
$$|E\cap B_\e|\ge{\half} |B_1| \e^n,\qquad 0<\e<\e_0.\eqdef{12w}$$  
It is no big surprise that even in the case of the Riesz potential $I_\a$ (i.e. $g\equiv1$) a family of functions satisfying \eqref{12x}, with $A_g=|B_1|$,  will be obtained by a  suitable modification of the usual extremal Adams family
$$\phi_\e(y)=\cases{|y|^{-\a} & if $\e\le |y|\le 1$\cr 0 & otherwise\cr},\qquad 0\le\e<1$$
It is clear that some modification is necessary, due to the integrability requirements on $(I_\a \phi_\e)^\na$ at infinity. Indeed,  $|x-y|^{\a-n}\sim |x|^{\a-n}$ when  $|x|$ is large, and this implies that for any $f$ compactly supported in $\Rn$,  $I_\a f$ cannot be in $L^\na(\R^n)$ for $n/2\le \alpha<n$,  if $\int_{\sR^n} f\neq0$. The same considerations can be made for general kernels $g$.

For slightly more clarity we will  prove the result in the scalar case first, but the modifications in the vector-valued case are simple, and will be indicated after the proof of the scalar case. So, assume that $g\in C^{2n}(\Sn)$ and for $|y|\le 1$  and $|x|\ge 3$ from  Taylor's formula, there exists $\theta=\theta(x,y)\in(0,1)$ such that 
$$\eqalign{g(x-y)&=|x|^{\a-n}\sum_{k=0}^{2n-1} {1\over k!} d^k g\Big(x^*,-{y\over|x|^2}\Big)+{|x|^{\a-n}\over (2n)!}d^{2n} g\Big(x^*- {\theta y\over|x|^2},-{y\over|x|^2}\Big)\cr&:=|x|^{\a-n}\sum_{k=0}^{2n} {1\over k!} p_k(x,y).\cr}\eqdef{g}$$
Now, $p_k(x,y)$ is a polynomial of order $k$ in the $y$ variable for $k\le 2n-1$, and $|p_k(x,y)|\le C|x|^{k}$, for  $0\le k\le 2n$, 
for some constant $C$ independent of $y$ and $x$ in the given range.

Next, let $\P_m$ be  the space of polynomials of degree $m$ in the unit ball $B_1$ of $\Rn$, a subspace of $L^2(B_1)$. Let $\{v_1,....,v_N\}$ be an orthonormal basis of $\P_m$, with $v_1=|B_1|^{-1/2}$. If $P_m$ denotes the projection of $L^2(B_1)$ onto $\P_m$, then $P_m$ has integral kernel
$$P_m(y,z)=\chi_{B_1}^{}(y)\sum_{k=1}^N v_k(y)v_k(z)$$
which is pointwise uniformly bounded on $B_1\times B_1$ (with bound depending on $m$). 

Let $$\phi_\e(y)=\cases{\big( g(y^*)\big)^{\a\over n-\a}|y|^{-\a} & if $\e\le |y|\le 1$\cr 0 & otherwise\cr},\qquad 0\le\e<1$$
and  consider the functions on the unit ball 
$$\wtilde \phi_\e=\phi_\e-P_{2n}\phi_\e,\eqdef{12g}$$
which are orthogonal to every polynomial of order up to (and including) $2n$. Since $|P_{2n}\phi_\e|\le C\|\phi_\e\|_1=C$, then for all $\e>0$ small enough
$$\|\wtilde \phi_\e\|_{n/a}^{\na}=\|\phi_\e\|_\na^\na+O(1)=nA_g\log{1\over \e}+O(1).\eqdef{12gg}$$
If $|x|\ge 3$, then 
$$T_g\wtilde \phi_\e(x)={|x|^{\a-n}\over (2n)!}\int_{B_1}\wtilde \phi_\e(y)\,p_{2n}(x,y)dy,$$
hence  $|T_g\wtilde \phi_\e(x)|\le C|x|^{\a-2n-1}$ and  
$$\int_{|x|\ge 3}|T_g\wtilde \phi_\e(x)|^{\na}\le C.\eqdef{12t}$$
To handle the case $|x|\le 3$ we note that $|\wtilde \phi_\e|\le |\phi_0|+C\chi_{B_1}^{}$, and $|T_g\wtilde \phi_\e|\le T_{|g|}|\phi_0|+C\in L^{\na}(B_3)$, which can be checked for example via the standard O'Neil inequality:
$$(T_{|g|}|\phi_0|)^{**}(t)\le \nsa A_g^{n-\a\over n}t^{{\a\over n}-1}\int_0^t\phi_0^*(u) du+A_g^{n-\a\over n}\int_t^{\infty}u^{\asn-1}\phi^*_0(u)du\le C\Big(1+\log^+{|B_1|\over t}\Big),$$
if $0<t\le 3^n|B_1|$. Alternatively, note that $\phi_0\in L^p(\Rn)$ for each $p<\nsa$  (and $p>1)$, so that $T_g \phi_0\in L^q$
with $q=p^{-1}-\a/n$, hence $T_g\phi_0\in L^{\na}(B_3)$ if we pick any $p>1$ with  ${n\over 2\a}<p<\nsa$.
 In summary, we have that
$$\|\wtilde \phi_\e\|_\na^\na+\|T_g\wtilde \phi_\e\|_\na^\na=nA_g\log{1\over\e}+O(1)$$
for all $\e>0$ small enough. Lastly, we estimate for $|x|\le \e/3$
$$|T_g\wtilde\phi_\e(x)|\ge T_g\phi_\e(x)-C\ge T_g\phi(0)-|T_g\phi_\e(x)-T_g\phi_\e(0)|-C\ge nA_g \log{1\over\e}-C,\eqdef{12ggg}$$
(using for example \eqref{22k} with $R=\e/3$.).
If we now define
$$\psi_\e={\wtilde\phi_\e\over\Big(\|\wtilde \phi_\e\|_\na^\na+\|T_g\wtilde \phi_\e\|_\na^\na\Big)^\an}$$
then $\|\psi_\e\|_\na^\na+\|T_g\psi_\e\|_\na^\na=1$, and 
$$|T_g\psi_\e(x)|^\nna\ge nA_g\log{1\over\e}-C\Big(nA_g\log{1\over\e}\Big)^{n-\a\over n},\qquad |x|\le {\e\over3}.$$
Therefore,
$$\eqalign{\int_{E\cap B_{\e/3}}\exp\bigg[&{1+\delta\over A_g}\,|T_g \psi_\epsilon(x)|^{n\over n-\a}\bigg]dx
\ge|E\cap B_{\e/3}|\exp\bigg[(1+\delta)n \log{1\over\e}-C\Big(\log{1\over\e}\Big)^{n-\a\over n}\bigg]\cr& \ge C\exp\bigg[\delta n\log{1\over \e}-C\Big(\log{1\over\e}\Big)^{n-\a\over n}\bigg]\to+\infty\cr}$$
which  proves \eqref{12x}, in the case $g$ scalar. In the vector-valued case the proof is completely similar. 
First write an expansion as in \eqref{g} where each $p_k$ is a vector-valued polynomial whose components correspond to the Taylor's formula of each $g_j(x-y)$. 
Then define
$$\phi_\e(y)=\cases{\big( g(y^*)|g(y^*)|^{{\a\over n-\a}-1}|y|^{-\a} & if $\e\le |y|\le 1$\cr 0 & otherwise\cr},\qquad 0\le\e<1$$
and $\wtilde\phi_\e=\phi_\e-P_{2n}\phi_\e$, where $P_{2n}$ acts component-wise. The rest of the argument is exactly as in the scalar case. This concludes the proof of Theorem 7.
\endpf
\bigskip
\noin{\bf Note.} We would like to emphasize that the precise steps where Ruf's condition is used are in \eqref{usaruf1} and \eqref{usaruf2}. Those steps also make it clear that the reason why the proof of Theorem 7 fails if one uses the condition $\|f\|_\na^{q\na}+\|T_gf\|_\na^{q\na}\le1$ is that the inequality $\|f\|_\na^{qn/\a}\ge\|f\chi_B^{}\|_\na^{qn/\a}+\|f\chi_{B^c}^{}\|_\na^{qn/\a}$ is only true for $q\le1$ (being trivially an equality when $q=1$).

\bigskip
 
\centerline{\bf A simple proof  of an Adams inequality \`a la Adachi-Tanaka}

\bigskip

In this section we prove a special case of Corollary 2, namely that for given $\theta\in(0,1)$ there is $C=C(\a,\theta,g,n)$ such that 
$$\int_E \exp\Big[{\theta\over A_g}|T_g f(x)|^\nna\Big]dx\le C\eqdef{AT1}$$
for all $f$ such that $\|f\|_\na\le1$ and $\|T_g f\|_\na\le 1$.
This proof does not require the full force of Theorem 7, but uses instead only (A) in Lemma 10 and Lemma 14. Suppose $\| f\|_\na\le 1$, $\|T_g f\|_\na\le 1$, then $|T_g (\theta^{n-\a\over n} f)|\le |T_g (\theta^{n-\a\over n} f_\ell)|+C$ by Lemma 14, hence \eqref{AT1} follows at once  from (A) of Lemma 10 applied to $\theta f_\ell$.

Obviously  \eqref{AT1}  follows under the more restrictive  condition $$\|f\|_\na^{q\na}+\|T_g f\|_\na^{q\na}\le~1,\qquad q\ge 1,$$ but we will now show that under this condition with $1<q\le+\infty$ the inequality fails if $\theta=1$. 
Consider the functions
$$\wtilde \phi_{\e,r}(x)=r^{-\a} \wtilde\phi_\e(x/r)$$
where the $\wtilde \phi_\e$ are defined as in \eqref{12g}. Then $\supp\, \wtilde \phi_{\e,r}=B_r\setminus B_{r\e},$ and $$ \|\wtilde\phi_{\e,r}\|_\na=\|\wtilde \phi_\e\|_\na,\quad T_g\wtilde\phi_{\e,r}(x)=(T_g\wtilde\phi_\e)(x/r),\quad \|T_g\wtilde \phi_{\e,r}\|_\na=r^{\a}\|T_g\wtilde \phi_{\e}\|_\na.$$
For the rest of this argument choose $\e$ so that
$$\log{1\over \e^n}=r^{nq'}\eqdef{eps}$$
which is possible since $q>1$. If we define
$$\psi_{\e,r}={\wtilde\phi_{\e,r}\over\Big(\|\wtilde \phi_{\e,r}\|_\na^{q\na}+\|T_g\wtilde \phi_{\e,r}\|_\na^{q\na}\Big)^{\a\over qn}}$$ then $\big(\| \psi_{\e,r}\|_\na^{q\na}+\|T_g \psi_{\e,r}\|_\na^{q\na}\big)^{1/q}=1$ and,  owing to \eqref{12gg}, \eqref{12ggg}, for all $r$ large and $|x|\le r\e/3$ we get
$$|T_g \psi_{\e,r}(x)|\ge{A_g\log\e^{-n}-C\over \Big[ \Big(A_g\log \e^{-n}\Big)^q+Cr^{qn}\Big]^{\a\over qn}}\ge \Big(A_g\log\e^{-n}\Big)^{n-\a\over n}\bigg(1-C\Big({r^n\over\log\e^{-n}}\Big)^q\bigg)$$
$$|T_g \psi_{\e,r}(x)|^\nna\ge A_g\log \e^{-n}-C{r^{nq}\over \big(\log \e^{-n}\big)^{q/q'}}= A_g\log \e^{-n}-C$$
and
$$\int_{E\cap B_{r\e/3}}\exp\bigg[{1\over A_g}\,|T_g \psi_{\epsilon,r}(x)|^{n\over n-\a}\bigg]dx\ge C r^n\exp\bigg[-C{r^{nq}\over \big(\log \e^{-n}\big)^{q/q'}}\bigg]= Cr^n\to+\infty. $$
\endpf
\bigskip

\centerline{\bf Proof of Corollary 8}\bigskip
Assume the Adachi-Tanaka condition \eqref{103}. Clearly it suffices to prove  estimate \eqref{102} for $0<\theta_0\le\theta<1$.
For any $\lambda>0$ if $f_\lambda(x)=\lambda^\a f(\lambda x)$ then using \eqref{D1} and \eqref{D2}
$$\int_E \exp\Big[{\theta\over A_g}|T_g f(x)|^\nna\Big]dx=\lambda^n\int_{E/\lambda} \exp\Big[{1\over A_g}|T_g (\theta^{n-\a\over n} f_\lambda)(x)|^\nna\Big]dx\eqdef{D3}$$
and  
$$\|\theta^{n-\a\over n} f_\lambda\|_\na^\na+\|T_g(\theta^{n-\a\over n} f_\lambda)\|_\na^\na=\theta^{n-\a\over\a}\Big(\|f\|_\na^\na+\lambda^{-n}\|T_g f\|_\na^\na\Big)\le\theta^{n-\a\over\a}\big(1+\lambda^{-nq'}\big)^{1/q'}=1$$
 for $\lambda=\lambda(\theta)=\big(\theta^{-q'{n-\a\over \a}}-1\big)^{-{1\over nq'}}\ge \lambda(\theta_0)=1$, if we choose $\theta_0=2^{-{n\over q'(n-\a)}}$.
We can then apply Theorem 7 to estimate  \eqref{D3} with $C\lambda^n\le C(1-\theta)^{-1/q'}$.

Regarding the sharpness statement,  the family $\{\psi_{\e,r}\}$ in the previous proof satisfies
$$\int_{E\cap B_{r\e/3}}\exp\bigg[{\theta\over A_g}\,|T_g \psi_{\epsilon,r}(x)|^{n\over n-\a}\bigg]dx\ge C r^n\e^{(1-\theta)n}=Cr^ne^{-(1-\theta)r^{nq'}},$$
therefore it is enough to choose $r=(1-\theta)^{-{1\over nq'}}$.\endpf

\bigskip
\centerline{\bf  Proof of Theorem 1.}\bigskip  The Adams inequality for $I_\a*f$ under Ruf's condition is a special case of Theorem 7.

Regarding (b),  is enough to prove the results for $u\in C_c^\infty(\Rn)$. For any $\a<n$, $\a$ even,  we can write
$$u(x)=c_\a I_\a* f(x)=\irn c_\a |x-y|^{\a-n} f(y)dy,\quad f=(-\Delta)^{\a\over2} u\in C_c^\infty(\Rn)$$
and  Ruf's condition \eqref{12b}  for $u$ translates directly into Ruf's condition \eqref{101} for $f=\Delta^{\a\over2}u$, so  part a)  applies, yielding inequality \eqref{12c}.

If $\a<n$ is an odd integer, writing $u=c_{\a+1}I_{\a+1}(\Delta^{\a+1\over2} u)$ and integrating by parts gives
$$u(x)=J_\a f(x)=\irn c_{\a+1}(n-\a-1)|x-y|^{\a-n-1}(x-y)\cdot f(y)dy,\quad f=\nabla(-\Delta)^{\a-1\over2}u.\eqdef{J}$$

The kernel of  $J_\a$, changes sign component-wise, however we will verify that the alternate pointwise condition \eqref{12q} of Theorem 7 holds, for our given $f$. Specifically, if $J_\a^+$ is the potential with kernel $c_{\a+1}(n-\a-1)|x-y|^{\a-n}$ and $f$ is as in \eqref{J}, 
then we can prove that for each $a\in \R^n$ 
$$|J_\a f(x)|=|u(x)|\le J_\a^+|f\chi_{|y-a|\le 2}^{} |(x)+C\|J_\a f\|_\na,\qquad|x-a|\le 1.\eqdef{Q}$$
It is enough to prove this for $a=0$ on the function $u_a(x)=u(x -a)$, so WLOG we can assume $a=0$.

Indeed, pick any smooth $\phi\in C_0^\infty(\Rn)$ such that $0\le|\phi|\le 1$ and 
$$\phi(y)=\cases {1 & if $|y|\le{3\over2}$\cr 0 & if $|y|\ge2$\cr}$$
and write, using Leibinz's rule and integration by parts (differentiations are in the $y$ variable),  
$$\eqalign{u(x)\phi(x)&=J_\a\big(\nabla\Delta^{\a-1\over2}(u\phi)\big)(x)=\int_{|y|\le2} c_{\a+1}\Delta|x-y|^{\a+1-n}\Delta^{\a-1\over2}(u\phi)(y)dy\cr&=
\int_{|y|\le2} c_{\a+1}\phi(y)\Delta|x-y|^{\a+1-n}\Delta^{\a-1\over2}u(y)dy+\cr&\hskip3em +
\int_{|y|\le2} c_{\a+1}\Delta|x-y|^{\a+1-n}\sum_{|k|+|h|\le {a-1}\atop |k|>0}b_{k,h,\a}(D^k\phi) (D^{h}u)\cr&
= -\int_{|y|\le2}  c_{\a+1}\phi(y)\nabla |x-y|^{\a+1-n}\cdot\nabla\Delta^{\a-1\over2}u(y)dy- \cr& \hskip3em 
-\int_{|y|\le2}  c_{\a+1}\Big(\nabla\phi(y)\cdot\nabla |x-y|^{\a+1-n}\Big)\Delta^{\a-1\over2}u(y)dy
+\cr&\hskip3em +\int_{|y|\le2} c_{\a+1}\sum_{|k|+|h|\le {a-1}\atop |k|>0}(-1)^{|h|}b_{k,h,\a}D_y^{h}\Big(\Delta|x-y|^{\a+1-n}D_y^k\phi(y)\Big)u(y)dy\cr}$$
 where $k=(k_1,...,k_n)$, $h=(h_1,...,h_n)$ are  multiindices, and the constants $b_{k,h,\a}$ are so that
$$\Delta^{{\a-1\over2}} (u\phi)=\sum_{|k|+|h|\le{\a-1}}b_{k,h,\a}(D^k\phi) (D^{h}u).$$
With further integrations by parts we can write
$$\eqalign{u(x)\phi(x)&
= -\int_{|y|\le2}  c_{\a+1}\phi(y)\nabla |x-y|^{\a+1-n}\cdot\nabla\Delta^{\a-1\over2}u(y)dy+ \cr& \hskip3em 
+ \sum_{0<|k|+|h|\le {\a+1}}c_{k,h,\a}\int_{|y|\le2} \Big(D_y^{h}|x-y|^{\a+1-n}D_y^k\phi(y)\Big)u(y)dy\cr}\eqdef{L} $$
for some other coefficients $c_{h,k,\a}$. Note that the derivatives of the function $\phi$ in the second term of \eqref{L}  all have positive order.
Now, for $|k|>0$ we have $\supp D^k\phi \subseteq \{{3\over2}\le |y|\le 2\}$, and for any fixed $x$ with $|x|\le 1$ the function  $y\to |x-y|^{\a+1-n}$ is $C^\infty$ outside the ball of radius ${3\over2}$, so that for all such $x$ 
$$\eqalign{|u(x)|&=|u(x)\phi(x)|\le\int_{|y|\le 2}c_{\a+1}\big|\nabla |x-y|^{\a+1-n}\big|\,\big|\nabla\Delta^{\a-1\over2}u(y)\big|dy+C\int_ {{3\over2}\le |y|\le2}|u(y)|dy\cr&\le \int_{|y|\le 2}c_{\a+1}(n-\a-1)| |x-y|^{\a-n}\big|\nabla\Delta^{\a-1\over2}u(y)\big|dy+C\|u\|_\na, \cr} $$
which is \eqref{Q}.

To prove that the constant in \eqref{101} is sharp in the case $\alpha<n$  with $\alpha$ even, it is enough to consider the functions $u_\e=c_\a I_\a*\psi_\e $, where  the $\psi_\e$ were constructed in the proof of the sharpness statement in Theorem 7. In the case $\alpha$ odd, we can take the same extremal family $\{u_\e\}\in W_0^{\a,\nsa}(B(0,1))$ used in the original proof by Adams (see also  [FM1] proof of Theorem 6). Essentially, if  $v_\e$   is a  smothing of the function
$$\cases{0 & if $|y|\ge{3\over4}$\cr \log{1\over|y|} & if $2\e\le |y|\le {1\over2}$\cr
\log{1\over\e} & if $|y|\le \e$\cr}$$
then it is easy to check that 
$$\|v_\e\|_\na\le C,\qquad \|\nabla\Delta^{\a-1\over2}v_\e\|_\na^\na=\omega_{n-1}^{-{n-\a\over \a}}\big((n-\a-1)c_{\a+1}\big)^{-\nsa}\log{1\over\e}+O(1),$$ 
and that the exponential integral in \eqref{12c} evaluated at the  functions $u_\e=v_\e(\|v_\e\|_\na^\na+\|\nabla\Delta^{\a-1\over2}v_\e\|_\na^\na)^{-\asn}$ can be made arbitrarily large if the exponential constant is larger than $\gamma_{\a,n}$.

\endpf
  
\smallskip
\centerline{\bf Proof of Corollary 2}
\bigskip
This proof is identical to the one of Corollary 8. Given a function $u$ satisfying $\|u\|_\na^{qn/\a}+\|\nabla^\a u\|_\na^{q\na}\le 1$, we consider the functions $u_\lambda(x)=u(\lambda x)$, for $\lambda>0$, which satisfy  
$$\|u\|_\na^\na=\lambda^{-n}\|u\|_\na^\na,\qquad \|\nabla^\a u_\lambda\|_\na=\|\nabla^\a u\|_\na$$
and we choose $\lambda=\lambda(\theta)$ as we did in the proof of Corollary 8 to obtain inequality \eqref{102}.

The proof of the sharpness statement is also similar. All we need to do is take  the family $\{u_\e\}$ that extremizes \eqref{101} (with $u_\e\in W_0^{\a,\nsa}(\Rn)$)  and consider the family $\{u_{\e,r}\}$, where $u_{\e,r}(x)=u_\e(rx)$, with $\log\e^{-n}=r^{nq'}$ and $r=(1-\theta)^{-{1\over nq'}}$.
\endpf
\smallskip
\centerline{\scaps 3. Proof of Theorem 3}

\bigskip
As we mentioned in the introduction, the proof of Theorem 3 is accomplished by making slight modifications to the proof in [FM1, Theorem 1], in order to take into account the integrability condition \eqref{108}, and by tracking down the various constants a little bit more carefully. For the convenience of the reader we will present here the beginning of the proof in enough details so that the role of \eqref{108} is highlighted, relegating the more technical part (Adams-Garsia's lemma) to the appendix.
\smallskip

First observe that \eqref{110} follows from \eqref{109} and Lemma 9. It is then enough  to prove that for each $f\in L^\bp(M)$ the function $Tf$ is well-defined, finite a.e., and satisfies 
$$\int_0^{\nu(N)} \exp\bigg[{\s\over A}|(Tf)^{**}(t)|^\beta\bigg]dt\le   Ce^{\s J}\big(1+J+\nu(N)\big),\eqdef{0z}$$
where $ C=C(\b,\b_0,\gamma,H,A,B)$, under the hypothesis \eqref{106}, \eqref{107}, \eqref{108}, and with $$\int_0^\infty (f^*)^{\b'}\le 1.\eqdef{0d}$$
Below, $C_j$ denotes a constant $\ge 1$, depending only on $A, B, \b,\s,p,H,\gamma$.
\smallskip
WLOG we can assume that $k$ and $f$ are nonnegative. Clearly $k_1^*(t)\le A^{\bi} (1+H)t^{-\bi}$ for $t>0$, so that  by the improved O'Neil inequality in [FM1, Lemma 2],  if $p$ is any fixed number such that
$$\max\bigg\{1,{\b(1-\s),\over \b-1}\bigg\}\le p<{\b\over\b-1}=\b'\eqdef{0a}$$
and
$${1\over q}={1\over\s\b}+{1\over\s}\bigg({1\over p}-1\bigg),\qquad q>p\eqdef{0b}$$
then there is $C_0$ such that for each $t>0$ 
$$(Tf)^{**}(t)\le C_0 t^{-{1\over q}}\int_0^{t^{1/\s}} f^*(u)u^{-1+{1\over p}}du+\int_{t^{1/\s}}^\infty k_1^*(u)f^*(u)du.\eqdef{1}$$

\noin{\bf Note.} The general results in Part I of [FM1] were proved under the assumption that $\mu(M)$ is finite. However such condition is not necessary for the validity of the improved O'Neil inequality [FM1, eq. (19)], and neither is the integrability condition \eqref{108}. 
\smallskip
\noin{\bf Note.} When $\s=1$ we can take $p=1$ and $q=\b$ in \eqref{1}.\smallskip

If $t\ge 1$ 
 then H\"older's inequality and \eqref{0d} imply
$$(Tf)^{**}(t)\le \bigg(\int_0^{t^{1/\s}} C_0^\b t^{-{\b\over q}}u^{\b\big({1\over p}-1\big) }du+\int_{t^{1/\s}}^\infty \big(k_1^*(u)\big)^\b du\bigg)^\bi\le\bigg( {C_0^\b\over\b\big({1\over p}-1\big)+1}+AJ\bigg)^\bi.$$
If instead $t<1$, then
$$\eqalign{(Tf)^{**}(t)&\le \bigg(\int_0^{t^{1/\s}} C_0^\b t^{-{\b\over q}}u^{\b\big({1\over p}-1\big) }du+\int_{t^{1/\s}}^1 \big(k_1^*(u)\big)^\b du+\int_{1}^\infty \big(k_1^*(u)\big)^\b du\bigg)^\bi\cr&\le\bigg({C_0^\b\over\b\big({1\over p}-1\big)+1}+A\int_{t^{1/\s}}^1\big(1+H(1+|\log u|)^{-\gamma}\big)^\b \,{du\over  u} +AJ\bigg)^\bi\cr&\le\bigg(C_1+A\log{1\over t^{1/\s}}+AJ\bigg)^{\bi}.\cr}$$
Hence, for any $t>0$
$$(Tf)^{**}(t)\le \bigg(C_1+A\log^+{1\over t^{1/\s}}+AJ\bigg)^{\bi}$$
which shows in particular that $Tf$ is finite a.e.  The same inequality also shows that 
$$\int_{1}^{\nu(N)} \exp\bigg[{\s\over A}|(Tf)^{**}(t)|^\beta\bigg]dt\le Ce^{{\s }J} \big(\nu(N)- 1\big)^+.$$
On the interval $[0,1]$ unfortunately this simple argument fails  and we need to refine the more sophisticated analysis   in [A] and [FM1].  Make the change $u=v^{1\over\s}$ in \eqref{1} and then the changes
$$v=e^{-x},\quad t=e^{-y},\quad \phi(x)=\bigg({1\over\s}\bigg)^\bpi f^*\big(e^{-{ x\over\s }}\big) e^{-{\b-1\over\s\b}x}\eqdef{1q}$$
to obtain 
$$\eqalign{&\int_0^{ 1} \exp\bigg[{\s\over A}|(Tf)^{**}(t)|^\beta\bigg]dt\le 
\int_0^{1} \exp\bigg[\bigg({1\over\s}\bigg)^{\b\over\b'}\bigg(C_2^\bi t^{-{1\over q}}\int_0^t f^*\big(v^{1\over \s}\big)v^{-1+{1\over\s p}}dv+\cr& +\int_t^{1} \big(1+H_1(1+|\log v|)^{-\gamma}\big)v^{{\b-1\over\s\b}-1} f^*\big(v^{1\over \s}\big)dv+A^{-\bi}\int_{1}^\infty k_1^*\big(v^{1\over \s}\big)f^*\big(v^{1\over \s}\big)v^{{1\over\s}-1}dv\bigg)^\b\,\bigg] dt\cr
&=\int_0^\infty e^{-F(y)}dy \cr
}\eqdef{1c}$$
(here $C_2=C_0^\b/A$, $H_1=H/\sigma$)
where for each fixed $y\ge 0$
$$F(y)=y-\bigg(\int_{-\infty}^\infty g(x,y)\phi(x)dx\bigg)^\b$$
$$g(x,y)=\cases{A^{-\bi}k_1^*(e^{-{x\over\s}})e^{-{x\over\s\b}} & if $x\le0$\cr
1+H_1(1+|x|)^{-\gamma} & if $0<x\le y$\cr
C_2^{\bi}e^{y-x\over q} & if $y<x$.\cr}$$
The next technical step is to run the Adams-Garsia machinery to prove that 
$$\int_0^\infty e^{-F(y)}dy\le C(1+J)e^{\sigma J}.$$
The remaining details  are given in the Appendix.\endpf

\medskip

\eject
\smallskip
\pf Proof of Corollary 4. In Theorem 3 let $N=E$ and $M=\R^n$ with the Lebesgue measure, and let $k(x,y)=K_\a(x-y)$, $\beta=\nna$. The proof that \eqref{111} implies \eqref{106} and  \eqref{107} for small $t$, and therefore for $t\in (0,1]$,  has been done in [FM1, Lemma 9]. Note that the proof there was done in the case $g$ bounded on the sphere, but it works even in our more general hypothesis.

It is enough to check that  \eqref{112} implies   \eqref{108}  (from which  \eqref{107} follows for all $t$).
The proof of this fact is  straightforward. Let $|K_\a(x)|\le M$ for $|x|\ge R$, and let 
$$\wtilde K_\a(x)=\cases{|K_\a(x)| & if $|x|\ge R$\cr M & if $|x|<R.$\cr}$$
If $\lambda(s)$ and $\wtilde\lambda(s)$ denote the distribution functions of $K_\a,\wtilde K_\a$ respectively, then $\wtilde\lambda(s)\ge \lambda(s)$ for $s< M$, and $\wtilde\lambda(M^-)\ge |B_1|R^n$. Hence, if $ k_1^*,\,\wtilde k_1^*$ denote the rearrangements of $ K_\a,\, \wtilde K_\a$ resp., then $\wtilde k_1^*(t)\ge k_1^*(t)$ for $t\ge |B_1|R^n$.
Obviously,  $\wtilde K_\a\in L^\nna(\Rn)$, so \eqref{108} follows if $|B_1|R^n\le1$. If $|B_1|R^n>1$ then \eqref{108} still follows since $k_1^*\in L^\nna([1,|B_1|R^n])$, since $k_1^*(1)<\infty$.
This proves  inequality \eqref{112a}, and therefore \eqref{113}.

The proof of the sharpness statement is the same as that of [FM1, Theorem 8].\endpf

\pf Proof of Corollary 6.
Take  $k(x,y)=g(x-y)$, $\b={n\over n-\alpha}$, so that under the assumption \eqref b we have (see also  [FM1, Lemma 9])
$$k_1^*(t)=A_g^\nna  t^{-\nna},\qquad k_2^*(t)\le Q^{n-\a\over\sigma n} t^{-{n-\a\over\sigma n}},\qquad t>0$$
If we  take $\b=\nna,$ $\b_0=\sigma\b$, $N=E\subseteq \Rn$, with $\nu(E)<\infty$ and $M=F\subseteq \Rn$, with $|F|<\infty$, then Theorem 3 gives the result in \eqref{18}, since $J=\log^+|F|$.

The proof of the sharpness result is the same as the one  of Theorem 5,   with the exception that this time $\nu(E\cap B_{\e/3})\ge c_1\e^{\sigma n}$ (using the notation in that proof, where we are taking $x_0=0$ and $\e<r_0$).

If $\nu$ is the Lebesgue measure, then we have for $|F|=1$
$$\int_E \exp\bigg[{1\over A_g}|T_g f(x)|^{n\over n-\alpha}\bigg]dx\le C(1+ |E|),\eqdef{cor1}$$
therefore, using the dilation $f_\lambda(x)=\lambda^\a f(\lambda x)$ we have $\supp f_\lambda\subseteq F/\lambda$, and formulas  \eqref{D1} and \eqref{D2} with  $\lambda=|F|^{1/n}$ and \eqref{cor1} give \eqref{18imp}. 
\endpf
\eject

\bigskip\centerline {\scaps 4. Further consequences of Theorem 3 and Theorem 7}
\bigskip\centerline{\bf Moser-Trudinger inequalities for elliptic operators with constant coefficients.}
\bigskip

In this section we give some applications of Theorem 7 and Theorem 3. Specifically, we will consider Moser-Trudinger inequalities for more general elliptic operators with constant coefficients, and obtain sharp inequalities for homogeneous operators from Theorem 7, and for  some non homogeneous differential  and pseudodifferential operators as a consequence of Corollary 4 (which is itself a consequence of Theorem 3).

Let  us consider an  elliptic differential operator of  order $\alpha<n$ with constant coefficients
$$Pu=\sum_{|k|\le\a}a_k D^k u$$
where $k=(k_1,...,k_n)$ denotes a nonnegative multiindex in $\Z^n$, acting, say, on $C^\infty_c(\Rn)$. We will let
$$p(\xi):=P(2\pi i \xi)=\sum_{|k|\le\a}a_k (2\pi i\xi)^k.$$
For simplicity here we only consider the case $a_k\in \R$, in which case $\a$ is even and ``elliptic"  means that the strictly homogeneous principal symbol of $P$ satisfies
$$p_\a(\xi):=P_\a(2\pi i\xi):=(2\pi)^\a(-1)^{\a/2}\sum_{|k|=\a}a_k \xi^k\ge c_0|\xi|^\a,\qquad \xi\in \R^n$$
for some $c_0>0$.

\proclaim Theorem 15. a) Let $P=P_\alpha$ be a homogeneous elliptic operator of order $\a<n$ with constant coefficients. Then there exists $C>0$ such that for every $u\in W^{\a,\nsa}(\Rn)$ with 
$$\|u\|_\na^\na+\|Pu\|_\na^\na\le 1\eqdef{K5}$$
 and for all measurable $E\subseteq \Rn$ with $|E|<\infty$ we have
$$\int_E\exp\bigg[{1\over A_P^{}}|u(x)|^\nna\bigg]dx\le C(1+|E|)\eqdef{K4}$$
and
$$\int_{\sR^n}\exp_{[\nsa-2]}\bigg[{1\over A_P^{}}|u(x)|^\nna\bigg]dx\le C\eqdef{K4a}$$
where 
$$A_P^{}={1\over n}\isn |{g_P^{}}^{}(\omega)|^\nna d\omega,\qquad{g_P^{}}^{}(x)=\irn {e^{-2\pi i x\cdot \xi}\over p_\a(\xi)}d\xi,\eqdef{K3}$$
 in the sense of distributions.
Moreover, the exponential constant $A_P^{-1}$ in \eqref{K4}, \eqref{K4a} is sharp.\smallskip
\noin  b) If $P$ is a non-homogeneous  elliptic differential operator  with constant coefficients of order $\a$, with ${n\over2}\le \a<n$ and  with $p(\xi)\neq 0$ for $\xi\neq 0$,  then \eqref{K4}  holds for all $u\in W^{\a,\nsa}(\Rn)$ such that $\|Pu\|_\na\le 1$, and \eqref{K4a} holds under the additional condition $\|u\|_\na\le c_1$, some fixed $c_1>0$. The exponential constants are sharp. \par

\def\k{{\kappa}}

\pf Proof. 
 It is well known that $P$ has a fundamental solution,  given by a function ${K_P^{}}^{}$  which is $C^\infty$ outside the origin, and which is  formally the inverse Fourier transform of $1/p( \xi)$. If  $p(\xi)\neq0$ for $\xi\neq 0$, a concrete formula for ${K_P^{}}^{}$ can be written for example as follows: 
$${K_P^{}}^{}(x)=\irn {\eta(\xi)\over p(\xi)}\,e^{-2\pi i x\cdot\xi}d\xi+{(-1)^\ell \over(2\pi|x|)^{2\ell}}\irn \Delta^\ell\bigg({1-\eta(\xi)\over p(\xi)}\bigg)\,e^ {-2\pi i x\cdot\xi}d\xi,\eqdef {K1}$$
for $x\neq0$, were $\eta$ is a smooth cutoff which is equal 1 for $|x|\le 1$ and 0 for $|x|\ge 2$, for any $\ell>{n-\a\over2}$.
Using this formula, or other standard methods, it is possible to see that 
$${K_P^{}}(x)={g_P^{}}^{}(x^*)|x|^{\a-n}+O(|x|^{\a-n+\e}),\qquad |x|\le1\eqdef{K2}$$
where $x^*/|x|$, and where ${g_P^{}}$ is given as in \eqref{K3}
 (see also [FM1], formulas (67), (69)).

The proof of the theorem is therefore accomplished by first proving a sharp Adams inequality for the convolution operator ${K_P^{}}*f$, where $f$ is compactly supported in $L^\nsa(\Rn)$, and then by applying it using the representation  $u={g_P^{}}*(Pu)$, where $u\in C_c^\infty(\Rn)$.

If $P$ is homogeneous of order $\a$, i.e. in the above notation $P=P_\a$, then ${K_P^{}}^{}$ given by \eqref{K3} is homogeneous of order $\a-n$ and ${K_P^{}}^{}(x)={g_P^{}}^{}(x^*)|x|^{\a-n}$.
We could therefore apply Theorem 7 to obtain the Adams inequality on $\Rn$ (with Ruf condition) for the convolution operator ${K_P^{}}^{}* f$ if either ${g_P^{}}^{}$ does not change sign  on $S^{n-1}$, or else if condition \eqref{12q} is satisfied. Even though it seems like a natural condition,  the nonnegativity of ${g_P^{}}$ does not seem easy to establish in general, but it is possible to show that \eqref{12q} is true. Indeed, we have that  for $f=Pu$
$$|{K_P^{}}*f(x)|=|u(x)|\le |{K_P^{}}|*|f\chi_{|y-a|\le 2}^{}|(x)+C\|u\|_\na,\qquad  |x-a|\le1$$
the proof of which is a repetition of the proof of \eqref{Q}, but using the operator $P$ instead of $\nabla\Delta^{\a+1\over2}$.
This establishes \eqref{K4} in case $P$ is homogeneous.

 If instead $P$ is not homogeneous, then we can resort to Corollary 4 provided condition \eqref{112} is verified, i.e. that at infinity ${K_P^{}}$ is bounded and in $L^{\nna}$. 

In the formula for ${K_P^{}}$ given in \eqref{K1} it is clear that the second term is a rapidly decreasing function, therefore the problem is to show integrability and boundedness of the first term. Boundedness follows from the Riemann-Lebesgue lemma. In general, it is hard to establish the precise behavior of the first term in \eqref{K1} at infinity, however, if $n\le 2\a$ we can use the  Hausdorff-Young inequality and prove that the Fourier transform of $\eta/p$ (i.e. the first term of \eqref{K1})  is in $L^{\nna}(\Rn)$, by showing  $\eta/p\in L^\nsa(\Rn)$ i.e. $1/p\in L^{\nsa}(B_1)$. 

\proclaim Lemma 16.  Let $p$ be a real-valued elliptic  polynomial or order $\alpha$, with $p(0)=0$  and $p(x)>0$ for $x\neq0$. Then $p(x)^{-1}\in L^{\nsa}(B_1)$ if and only if $p$ is not homogeneous.
\par

\pf Proof of Lemma 16. Obviously, if $p$ is homogeneous of order $\a$ then $1/p$ cannot be in $L^\nsa(B_1)$. Suppose $p$ is not homogeneous and let $p_\a,p_{\k}$ be the highest and lowest order homogeneous parts of $p$, of orders $\alpha$ and $\k<\a$  respectively. Then we can write   $p=p_\a+q+p_\k$ and the hypothesis imply that for all $x\in\Rn$

$$ p_\a(x)\ge c_0|x|^\a,\qquad p_\k(x)\ge0,\qquad p(x)\ge c_1|x|^\a$$ where $c_0,c_1>0$.\def\o{\omega} 

(For the last inequality: first find $R>0$ so that the inequality is true for $|x|\ge R$. Next, the function 
$$\wtilde p(r,\o,z)={p(r\o)\over r^\kappa}-z r^{\a-\kappa}$$
 is continuous in $[0,R]\times S^{n-1}\times \R$, and strictly positive on $[0,R]\times S^{n-1}\times\{0\}$.
For each of $(r,\o)\in [0,R]\times S^{n-1}$ there is an open neighborhood $U_{r,\o}$ of $(r,\o)$ and a $\delta_{r,\o}>0$ such that $\wtilde p (r,\o,z)>0$ on $U_{r,\o}\times (-\delta_{r,\o}, \delta_{r,\o})$. Compactness finishes the proof.)

Note also that $p_\kappa(\o)=0$ on a set of zero measure on $S^{n-1}$.

Write
$$\int_{|x|\le1}{1\over p(x)^\nsa} dx=\isn d\o\int_0^1{r^{n-1}\over p\big(r\o\big)^\nsa}dr=\isn d\o\int_0^{1/p_\k(\o)}{p_\k(\o)^nr^{n-1}\over p\big(rp_\k(\o)\o\big)^\nsa}dr.$$
To ease a bit the notation assume that $\o\in S^{n-1}$ is fixed and writing $p_\a=p_\a(\o),\,p_\k=p_\k(\o)$ we have
$${p_\k^nr^{n-1}\over p\big(rp_\k\o\big)^\nsa}={p_\k^n r^{n-1}\over \big(r^\a p_\k^\a p_\a+q(rp_\k\o)+r^\k p_\k^{\k+1}\big)^\nsa}={p_\k^{\nsa(\a-\k-1)} r^{n-1-{n\k\over\a}}\over \big(r^{\a-\k} p_\k^{\a-\k-1} p_\a+r^{-\k}p_\k^{-\k-1}q(rp_\k\o)+1\big)^\nsa}.
$$
We can now choose $r_0>0$ such that for all $r\le r_0$ and all $\o\in S^{n-1}$
$$r^{\a-\k} p_\k^{\a-\k-1} p_\a+r^{-\k}p_\k^{-\k-1}q(rp_\k\o)+1\ge {1\over2}$$
(recall that $\k+1\le\a$ and the lowest homogeneous part of $q$ has order greater than $\k$).
Hence we can write 
$$\int_0^{1/p_\k}{p_\k^nr^{n-1}\over p\big(rp_\k(\o)\o\big)^\nsa}dr\le2^{-\nsa}\int_0^{r_0}p_\k^{\nsa(\a-\k-1)} r^{n-1-{n\k\over\a}}dr+\int_{r_0}^{1/p_\k} {p_\k^n r^{n-1}\over \big(c_1r^\a p_\k^\a\big)^\nsa}
\le C(1+\big|\log p_\k\big|\big)$$
Now,  the function $\log p_\k(\o)$ is integrable on the sphere. By homogeneity it is easy to check that this is equivalent to the local integrability of $\log p_\k(x)$, which follows from this general lemma:

\proclaim Lemma 17. If $p$ is any complex-valued polynomial in $\Rn$ then the function $\log|p|$ is locally integrable in $\Rn$.\par

We have not seen this result in the literature, so we will give here  a short proof.

\pf Proof of Lemma 17. Suppose $p$ has degree $m$. By a linear transformation $x_1=x_1'$, $x_j=x_j'+\lambda_j x_1'$, $j=2,3..,n$,  we can assume that $p(x)= x_1^m+a_{m-1} x_1^{m-1}+...+a_1 x_1+a_0$, where the $a_j$ are polynomials in $x_2,...,x_n.$ If $Q=[a,b]\times Q'$ is any cube in $\Rn$, then for fixed $x_2,...,x_n\in Q'$  the polynomial $p(x)$ has $m$ complex roots $\rho_k=\rho_k(x_2,...,x_n)$, $k=1,...,m$, which are all contained inside a fixed ball of radius $R$. Then the result follows from Fubini's theorem, since
$$\int_a^b \big|\log|t-\rho|\,\big| dt\le C(a,b,R),\qquad \rho\in\C,\; |\rho|<R.$$

Back to the proof of $1/p\in L^\nsa(\Rn)$, if $p$ has real coefficients then we can apply directly Lemma 16 to conclude. If $p$ has some complex-valued coefficients, then 
apply the lemma to the polynomial $|p(x)|^2$ which satisfies  the conditions of Lemma 16 with $\alpha$ replaced by $2\alpha$.

To prove the sharpness of the exponential constants, it is enough to take the family of functions $u_\e={K_P^{}}*\psi_\e$, where the $\psi_\e$ were constructed in the proof of Theorem 7.

\endpf
\bigskip

The restriction $n/2\le\a<n$  in b) is only needed in order to apply the Hausdorff-Young inequality, and could perhaps be removed with a more refined analysis of the kernel ${K_P^{}}$, that is the behavior at infinity of the first term in \eqref{K1}. For special non-homogeneous elliptic operators of any even order $\a<n$, it is possible to obtain enough information about the decay at infinity of $g_P^{}$. For example if $P=P_\a+Q_\a'$, where $0\le\a'<\a<n$ and $Q_\a'$ is elliptic and not identically 0, then it is possible to see that $K_p^{}(x)$ (that is the first term of \eqref{K1}) decays at least as $|x|^{\a'-n}$, which is enough to guarantee that $K_P\in L^\nna$ at infinity.
\bigskip
\vskip1em

\centerline{\bf  Moser-Trudinger inequalities in hyperbolic space}\bigskip

In this section  we obtain the sharp Moser-Trudinger inequalities for the higher order gradients on the hyperbolic space $\Hn$, as a  consequence of Theorem 3. Below, $\Hn$ will denote the hyperbolic space 
modeled by the forward sheet of the hyperboloid
$$\Hn=\{(x_0,x_1,...,x_n)\in\R^{n+1}:\, x_0^2-x_1^2-....-x_n^2=1,\; x_0>0\},$$
 endowed with the metric induced by the form
$$[x,y]=x_0y_0-x_1y_1-...-x_ny_n$$
and with distance function $d(x,y)={\rm {arccosh}} [x,y]$. One can introduce polar coordinates on $\Hn$
via 
$$x=(\cosh r,\sinh r\, \xi),\qquad r\ge0,\,\xi\in S^{n-1}$$
and in these coordinates the metric and the volume element are  given as 
$$ds^2=dr^2+\sinh^2 r d\xi^2,\qquad d\nu(x)=(\sinh r)^{n-1}drd\xi.$$
The Laplace-Beltrami operator on $\Hn$ is denoted as $\Delta_{\sHn}$, and in polar coordinates is given as 
$$\Delta_{\ssHn}={\p^2\over \p r^2}+(n-1)\coth r{\p\over\p r}+{1\over \sinh^2 r}\Delta_{S^{n-1}},$$
whereas the gradient $\nabla_\ssHn$ is computed as 
$$\nabla_\ssHn={\partial\over\partial r}+{1\over \sinh^2 r}\nabla_{S^{n-1}}.$$

The Sobolev space $W^{\a,p}(\Hn)$ of integer order $\alpha$ is  defined in the standard way via the covariant derivatives $\nabla^k$: it is the closure of the space of $C^\infty$ functions $\phi$ such that
$$\|\phi\|_{\a,p}:=\sum_{k=0}^\alpha \|\nabla^k\phi\|_p<\infty$$
where $\|\cdot\|_p$ denotes the norm in $L^p(\Hn,\nu)$.
As it turns out, on $\Hn$ it is enough to use the highest order derivatives in order to characterize the Sobolev space. In particular, if we define the higher order gradient on $\Hn$ as
$$\nabla_{\ssHn}^\a=\cases{\nabla_{\ssHn}(-\Delta_{\ssHn})^{\a-1\over2} & if $\alpha$ odd\cr (-\Delta_{\ssHn})^{\alpha\over2} & if $\alpha$  even,\cr}$$
then one has that $\|\nabla_\ssHn^\alpha u\|_p$ is an equivalent norm on $W^{\alpha,p}(\Hn)$. In particular, note that we have the Poincar\'e-Sobolev inequality 
$$\|u\|_p\le C\|\nabla_\ssHn^\alpha u\|_p,\qquad u\in W^{\a,p}(\Hn).$$
This inequality is proved in [Mancini-Sandeep-Tintarev] in the case of the gradient in the ball model (really a consequence of Hardy's inequality) and for even $\alpha$ in [Tat].

In this setup  sharp versions of the Moser-Trudinger inequality for $W^{\alpha,\nsa}(\Hn)$ are only known in the case $\alpha=1$ for the gradient ([MS], [MST], [LT1], [LT2]), and with the same sharp constant as in the Euclidean case. In the following theorem we give the general version of this result for arbitrary $\alpha:$
\medskip\
\proclaim Theorem 18. For any integer $\alpha$  with  $0<\alpha<n$ there exists a constant $C=C(\alpha,n)$ such that for every  $u\in W^{\a,\nsa}(\Hn)$ with 
$\|\nabla_\ssHn^\alpha u\|_\na\le 1$, and for all measurable $E$ with $0<\nu(E)<\infty$ we have 
$$\int_E\exp\Big[ \gamma_{n,\a}|u(x)|^{n\over n-\a}\Big]d\nu(x)\le C\big((1+\nu(E)\big),\eqdef {Z2}$$
and
$$\int_\sHn\exp_{[\nsa-2]}\Big[ \gamma_{n,\a}|u(x)|^{n\over n-\a}\Big]d\nu(x)\le C,\eqdef {Z2'}$$
  and the constant $\gamma_{n,\a}$ is sharp.

\par

\pf Proof. If $\alpha$ is even, the operator $(-\Delta_\ssHn)^{\a\over2}$ has a fundamental solution given by a  kernel of type $H_\a\big(d(x,y)\big)$, where $H_\a$ is positive and satisfies
$$H_\a(\rho)=c_\a \rho^{\a-n}+O(\rho^{\a-n+\e}),\qquad  \rho<1\eqdef{Z3}$$
(with the same $c_\a$ as in the Euclidean Riesz potential), and  $$H_\a(\rho)\le c_\a'\rho^{-1+{\a\over2}}e^{-(n-1)\rho}, \qquad \rho\ge1,\eqdef{Z4}$$
some $c_\a'>0$. These asymptotic estimates follows in a straightforward manner from the known formula for the fundamental solution of the Laplacian (see for example [CK])
$$H_2(\rho)={1\over\omega_{n-1}} \int_\rho^\infty {dr\over (\sinh r)^{n-1}}$$
using iterated integrations and the known addition formulas for the Riesz potential on $\R^n$. (In [BGS]  asymptotic formulas are derived  for general $\a$, using the  Fourier transform.)

It is now easy to check that \eqref{Z3} implies that in the measure space $(\Hn,\nu)$ we have  $$H_\a^*(t)=c_\a|B_1|^{n-\a\over n}t^{-{n-\a\over n}}+O(t^{-{n-\a\over n}+\e}),\quad t\le 1$$
while \eqref{Z4} implies that $H_\a\big(d(\cdot,O)\big)\in L^{\nna}\cap L^\infty\big(\{x:d(x,O)\ge1\}, \nu\big)$ (where $O=(1,0)$) and hence
$$\int_1^\infty (H_\a^*)^{n\over n-\a}dt<\infty.$$
Thus, we are in a position to apply Theorem 3 in  order to obtain \eqref{Z2} for $\alpha$ even, simply by writing $u(x)=\int H_\a(d(x,y))f(y)d\nu(y)$, with $f= \nabla_\ssHn^\a u$, for any $u\in C_0^\infty(\Hn)$.

If $\alpha$ is an odd integer, then we write 
$$u(x)=\int_\sHn \nabla_\ssHn H_{\a+1}(d(x,y))\cdot f(y)d\nu(y),\qquad f=\nabla_\ssHn^\a u$$
and use asymptotic estimates for $|\nabla_\ssHn H_{\a+1}|$,  which turn out to be the same exact estimates as in \eqref{Z3}, \eqref{Z4}, with $(n-\a-1)c_{\a+1}$ instead of $c_\a$.

The proof of the sharpness statement is identical to the one in the Euclidean case, namely we  let $v_\e$   to be  a   smoothing of the radial function
$$\cases{0 & if $r\ge{3\over4}$\cr \log{1\over r} & if $2\e\le r\le {1\over2}$\cr
\log{1\over\e} & if $r\le \e$.\cr}$$
Using local calculations as in [F, Prop. 3.6]  it is a routine task to check that 
if $\alpha$ is even then
 $$ \|\nabla_\ssHn^\a v_\e\|_\na^\na=\omega_{n-1}^{-{n-\a\over \a}}c_{\a}^{-\nsa}\log{1\over\e}+O(1),  $$
whereas is $\alpha$ is odd then the same estimate holds with $(n-\a-1)c_{\a+1}$ in place of $c_\alpha$.
From this estimate it is then clear that the exponential integral evaluated at the functions $u_\e=v_\e/\|\nabla_\ssHn^\a v_\e\|_\na$ can be made arbitrarily large  if the exponential constant is larger than $ \gamma_{n,\a}$. Note also that $\|v_\e\|_\na\le C$, so that $\|u_\e\|_\na\to0$ with $\e$, and the sharpness statement for the regularized inequality on $\Hn$ follows as well. \endpf\vskip1em

\centerline{\scaps 5. Further results and extensions} 
\bigskip
\centerline{\bf  Inequalities in critical potential spaces}\bigskip

The classical Sobolev embedding theorem states that  for $\alpha$ an integer between 0 and $n$ the Sobolev space $W^{\a,p}(\Rn)$ is embedded continuously into  $L^{np\over n-\a p}(\Rn)$, for $p<{\nsa}$, or, equivalently, that  the Riesz potential
 $I_\a$ maps $L^p(\Rn)$ continuously onto $L^{np\over n-\a p}(\Rn)$. At the critical index $p=\nsa$ the Moser-Trudinger inequality of Theorem 1 is a statement about the space $W^{\a,\nsa}(\Rn)$ being embedded   in an exponential class (sharply, in a suitable sense). Therefore it makes sense to ask whether there is a full analogue of this result at the level of Riesz potentials, namely whether the Adams inequality of Theorem 1 can be extended to a wider space of functions of $L^\nsa$, other than those having compact support, for any value of $\alpha\in (0,n)$.

Recall that for any $\a\in (0,n)$ the Riesz potential is well defined on the space 
$$L_c^{\na}(\R^n)=\{f\in L^\nsa(\Rn),\; {\rm supp } f\; {\rm compact}\},$$
and it is in $L^{\nsa}$  if $\alpha<n/2$. If $\alpha \ge n/2$ and $\phi\in L_c^{\na}(\R^n)$, then $I_\a*\phi$ is in $L^\nsa$ if and only if $\phi$ is orthogonal to the  Taylor polynomials of $|x^*-y|^{\a-n}$ (in the $y$ variable) up to order $\lfl 2\a-n\rfl$,  for all $x^*\in S^{n-1}$.  In this context it is then natural to   consider $I_\a$ as defined on the  linear subspace of $L^\nsa(\Rn)$
$$D(I_\a):=\{f\in L_c^{\na}(\R^n):\, I_\a*f\in L^\nsa(\Rn)\}.$$
The first observation is that $I_\a$ is not closed on $D(I_\a)$, not even if one defines $D(I_\a)$ using $C_c^\infty(\Rn)$ (or even $\S$, the space of Schwarz functions), rather than $L_c^\na(\Rn)$. Indeed, it's easy to exhibit a sequence $f_n\in C_c^\infty(\Rn)$ with $f_n\to f$ in $L^\nsa$, with $f \in \S\setminus C_c^\infty$, and $I_\a*f_n\to I_\a*f$.

 Hence, it is natural to ask whether or not $I_\a$ has a 
 smallest closed extension (in other words whether or not it is closable) as a densely defined operator from $L^\nsa$ to $L^\nsa$. We will now show that this is indeed the case, thereby allowing us to extend the notion of Riesz potential to a wider space than $D(I_\a)$, which we will call $\D_\a$. The space of  corresponding potentials $\U_\a:=I_\a(\D_\a)$, which we can call {\it critical Riesz potential space}, is the natural choice if one wishes to consider Moser-Trudinger inequalities  for abitrary powers of $-\Delta$. As it turns out, $\U_\a$ is a Banach subspace of the classical Bessel potential space $H^{\a,\nsa}(\Rn)$. We believe that these spaces actually coincide for all $\alpha<n$, but so far we can only prove it  for $\alpha$ even or $\alpha<n/2$.

We will now give a brief proof of the fact that $I_\a$ is closable,  since we were not able to find any references to this result.  We will in fact consider the more general homogeneous potential $T_g$ introduced at the beginning of section 2. We will treat the scalar case for simplicity, however the vector-valued case is treated similarly. 

For the rest of this section we will assume that $g:\Rn\setminus\{0\}\to \R$ is homogeneous of order $\alpha-n$ ($0<\alpha<n$), $g$ Lipschitz on $S^{n-1}$, and $T_gf=g*f$, which is well-defined in the space 
$$D(T_g):=\{f\in L_c^{\na}(\R^n):\, T_gf \in L^\nsa(\Rn)\}.$$

\proclaim Lemma 19. If $\{f_k\}\subseteq D(T_g)$ is such that $f_k\buildrel{ L^{\nsa}}\over \longrightarrow 0$ and $T_g f_k\buildrel{ L^{\nsa}}\over \longrightarrow h$, then $h=0$ a.e.

\par 
\pf Proof. If $h$ is not zero on a set of positive measure, we can assume that $\irn|h|^\nsa=1$
and that 
$$\int_{|x|\le R} |h|^\nsa\ge {3\over4},\qquad \int_{|x|\ge S} |h|^\nsa\le \e$$
for some $S>R>0$ and $\e$ small. Now consider
$$\varphi(x)={\rm sgn}(h)|h|^{n-\a\over\a}\chi_{B_R}^{}(x)$$
which is clearly in $L^{\nna}$ with $\|\varphi\|_\nna\le\|h\|_\na=1$, and has compact support, but its potential is not necessarily in $L^\nna$. For this to happen it is sufficient to normalize  $\varphi$ so that its mean  is zero, but we need to do this in a different way than the one used in the  proof of the sharpness statement of Theorem 1, which was localized inside a ball. 
\def\wvphi{{\widetilde\varphi}}

We let 
$$\wvphi(x)=\varphi(x)- \varphi(x-2S e_1)$$
and using \eqref{22ab}  we see that for $|x|\ge 4S$
$$\eqalign{|T_g\wvphi(x)|&\le \int_{|y|\le R}|g(x-y)-g(x-y+2Se_1)|\,| \varphi(y)|dy\cr&\le C_{\a,n}(2S)|x|^{\a-n-1}\int_{|y|\le R}|\varphi(y)|dy\le C|x|^{\a-n-1} \cr}$$
where $C$ depends on $R,S,\a,n$. Hence $T_g\wvphi$ is in $L^\nna$ for large $x$, and clearly this is also the case for small $x$.

Now we can say that $\wvphi\in L^\nna(\Rn)$ and $T_g\wvphi\in L^\nna(\Rn)$, and write
$$\eqalign{\irn& \wvphi h=\int_{B_R} |h|^\nsa-\int_{B_R}\varphi(x)h(x+2Se_1)dx\cr& 
\ge {3\over4}-\|\varphi\|_\nna\bigg(\int_{B_R}|h(x+2Se_1)|^\nsa dx\bigg)^\asn\ge {3\over4}-\e^{\asn}>{1\over2}\cr
}\eqdef{D1}$$
for $\e$ chosen small enough.

On the other hand, since $\wvphi,T_g\wvphi\in L^\nna$ and $f_k\to0,\,T_g f_k\to h$ in $L^\nsa$ we have 
$$\irn \wvphi h=\irn \wvphi \lim_k T_g f_k=\lim_k \irn \wvphi  T_g f_k=\lim_k \irn (T_g \wvphi)  f_k=\irn (T_g\wvphi)   \lim_k f_k=0$$
which contradicts \eqref{D1}\endpf


At this point we are in a position to apply a standard construction in order to close the operator $T_g$ see for ex. [Yosida Ch. II, Sect. 6]. Define
$$\D_g:=\big\{f\in L^\nsa(\Rn): \,\exists \{f_k\}\subseteq D(T_g),\,\exists h\in L^\nsa(\Rn)\; {\rm with}  \;f_k\buildrel{ L^{\nsa}}\over \longrightarrow f,\;T_g f_k\buildrel{ L^{\nsa}}\over \longrightarrow h\big\}\eqdef{D}$$
and because of Lemma 19 the function $h$ appearing in \eqref{D} is independent of the sequence $f_k$, and the potential $T_gf$ is well defined for $f$ in  $\D_g$, by letting $T_gf=h$. The operator thus defined is the smallest closed extension of $T_g$ as defined on $D(T_g)$, and the class $\D_g$ is the closure of $D(T_g)$ under the graph norm \eqref{D3}.

If the operator $T_g$ is injective, we define the critical potential space for the convolution operator $T_gf=g*f$ as 
$$\U_g=\{u\in L^\nsa(\Rn):\, u=T_g f,\; f\in \D_g\}$$
which is a Banach subspace of $L^\nsa(\Rn)$  endowed with the norm $\big(\|T_g^{-1}u \|_\na^\na+\|u\|_\na^\na\big)^{\a/n}$, or any other equivalent norm, such as  $\big(\|T_g^{-1}u \|_\na^{qn/a}+\|u\|_\na^{qn/\a}\big)^{\a/ qn}$, $q\in[1,\infty]$. The operator $T_g$ is obviously a continuous bijection between $\D_g$ and $\U_g$. If $g$ is a smooth function on the sphere (or differentiable enough times)
then the distributional F.T. $\widehat g$ is a homogeneous function of order $-\a$ and smooth on $\Rn\setminus\{0\}$; clearly $T_g$ is injective if $\widehat g\neq 0$ a.e. on $S^{n-1}$, for example like in the case of the Riesz potential.

Also, it is straightforward to check that $C_c^\infty(\Rn)\cap \D_g$ and the space of Schwarz functions in $\D_g$ are both dense in $\D_g$.

 In the special case $g(x)=c_\a|x|^{\a-n}$   we have that $T_gf=c_\a I_\a*f$ is the normalized Riesz potential, and the above procedure defines its smallest closed extension from a Banach space which we denote $\D_\a$, bijectively onto  the space 
$$\U_\a:=\{u\in L^\nsa(\Rn):\;u=c_\a I_\a* f,\, f\in \D_\a\},$$
 the corresponding   critical Riesz potential space. The inverse of $c_\a I_\a$, as defined on $\D_\a$,  is the fractional Laplacian  $(-\Delta)^{\alpha\over2}$, the sense that $u\in \U_\a$ if and only if $(-\Delta)^{\a\over2} u$ (which is defined distributionally as the inverse Fourier transform of $(2\pi|x|)^\a \widehat u$) is a function and it belongs to $L^\nsa(\Rn)$, and moreover $u=c_\a I_\a*\big((-\Delta)^{\a\over2} u\big)$. This can be easily verified using  density of $C_c^\infty(\Rn)\cap\D_\a$ in $\D_\a$, and the fact that in that space the statement is true.
  
Hence $\U_\a$ is a Banach subspace of the space 
$$ H^{\a,\nsa}(\Rn):=\{u\in L^\nsa(\Rn):\, (-\Delta)^{\a\over2}u\in L^\nsa(\Rn)\}$$
endowed with the norm $(\|u\|_\na^\na+\|(-\Delta)^{a\over2}u\|_\na^\na)^{\a/n}$. This space coincides with 
the classical Bessel potential space $\{u\in L^\nsa(\Rn):\, (I-\Delta)^{\a\over2}u\in L^\nsa(\Rn)\}$
with  norm \break $\|(I-\Delta)^{\a\over2} u\|_\na$,
where $(I-\Delta)^{\a\over2}u$ is the distributional inverse F.T. of $(1+|\xi|^2)^{\a\over2}\widehat u$.

 This  fact  is a consequence of the identities
$$(2\pi|x|)^{\a}=(1+\what h_1(x))(1+4\pi^2|x|^2)^{\a\over2},\;\;(1+4\pi^2|x|^2)^{\a\over2}=(1+\what h_2(x))\big(1+(2\pi|x|)^\a\big)$$
valid for some integrable functions $h_1,h_2$ (see [S, pp. 133-134]) which also imply that for all $f\in C_c^\infty(\Rn)$
$$c_1 (\|u\|_\na^\na+\|(-\Delta)^{a\over2}u\|_\na^\na)^{\a/n}\le \|(I-\Delta)^{\a\over2}u\|_\na\le c_2  (\|u\|_\na^\na+\|(-\Delta)^{a\over2}u\|_\na^\na)^{\a/n}.$$
It is well known that for $\alpha$ an integer  then $H^{\a,\nsa}=W^{\a,\nsa}$, the classical Sobolev space.

Regarding the connection between $\U_\a$ and classical Sobolev and Bessel potential spaces. When $\alpha$ is an even integer, it is  easy to see that $\U_\a=W^{\a,\nsa}$, for if $u\in W^{\a,\nsa}$ take $u_k\in C_c^\infty(\Rn)$ with $u_k\to u$ and $f_k:=(-\Delta)^{\a \over2} u_k\to f:=(-\Delta)^{\a\over2}  u$ in $L^\nsa$, and therefore, since $f_k$ also has compact support, we have $c_\a I_\a*f_k\to u$ and $u=c_\a I_\a * f$. It is also easy to handle the case $0<\a<n/2$:

\proclaim Proposition 20. For $0<\a<{n\over2}$ we have $\U_\a=H^{\a,\nsa}(\Rn)$.\par

\pf Proof. Since $C_c^\infty(\Rn)$ is dense in $H$ we only need to prove that if $u\in C_c^\infty(\Rn)$ then there exists $\{f_k\}\in L_c^\nsa(\Rn)$ with 
$f_k\to f$ and $T_\a f_k=I_\a*f_k\to u$, in $L^\nsa$, for some $f\in L^\nsa$. 

Let $u\in C_c^\infty(\Rn)$  then for all $x$ sufficiently large
$$|\Delta^{\a\over2} u(x)|\le C|x|^{-\a-n}.$$
For even $n$ this can be seen for example by writing $\Delta^{\a\over2}u=I_{n-\a}*\Delta^{n\over2} u$ and integrating by parts. For $n$ odd the argument is similar, starting from $\Delta^{\a\over2}u=\nabla I_{n-\a+1}*\nabla\Delta^{n-1\over2} u$.
Hence, for any $p\ge 1$ we have that $f:=\Delta^{\a\over2}u\in L^p(\Rn)$, and if 
$f_R=\chi_{B(0,R)}^{} \Delta^{\a\over2} u,$
then $f_R\to f$ in $L^p$, as $R\to+\infty$. The result follows since  $I_\a$ is continuous as an operator from  $L^{n\over 2\a}$ to $L^\nsa$.\endpf

In order to settle the identity $\U_\a=H^{\a,\nsa}(\Rn)$ for $\alpha\ge n/2$ a much more sophisticated argument than the one just presented seems to be needed.

We now observe that essentially every  theorem in this  paper can  formulated in terms of the spaces $\D_g$ and $\U_g$, (or $\D_\a$ and $\U_\a$  for arbitrary values of $\alpha$), by a simple limiting procedure. For example, we have the following:

\medskip\proclaim Theorem 21. \item{ a)}  For $0<\a<n$ the sharp Adams inequalities for the Riesz potential  \eqref{101} and \eqref{101b} hold for all $f\in \D_\a$ under the Ruf condition \eqref{100}. Likewise, if $g\in C^{2n}(S^{n-1})$ then the sharp  Adams inequality for the  general homogeneous potential \eqref{12} and \eqref{12aa} hold for all $f$ in the space $\D_g$, under the Ruf condition \eqref{11}. 
\item{b)} For $0<\a<n$ there is $C=C(\a,n)$ such that for all $u\in\U_\a$ with
$$\|u\|_\na^\na+\|(-\Delta)^{\a\over2}u\|_\na^\na\le1$$
we have
$$\int_{\sR^n}\exp_{[\nsa-2]}\bigg[{c_\a^{-\nna}\over|B_1|}|u(x)|^{n\over n-\a}\bigg]dx\le C$$
and the exponential constant is sharp.\par

\pf Proof. For a), it is enough to prove the inequality for $f$ having  Ruf norm equal to 1. Take a sequence $f_k$ of compactly supported functions such that $f_k\to f$ in the Ruf norm, and so that $T_g f_k\to T_g f$ a.e.. Apply the Adams inequality of Theorem 1 to  the normalized $f_k$, which also converge to $f$ in the Ruf norm, and use  Fatou's Lemma to conclude the proof. The proof of b) is immediate since $I_\a$ is a bijection between $\D_\a$ and $\U_\a.$\endpf

\smallskip
For the  case $\alpha=\half$, $n=1$ see also the recent paper  by Iula-Maalaoui-Martinazzi [IMM, Thm 1.5], which was proved by adapting Ruf's original argument. 
\bigskip\eject
\centerline{\bf  Inequalities for more general Borel measures}

\bigskip
The methods presented thus far allows us to obtain versions of the sharp inequalities in this paper when the non-regularized exponential is integrated against a positive Borel measure $\nu$ such that     
$$\nu\big(B(x,r)\big)\le Qr^{\sigma n},\qquad \forall x\in \Rn,\,\forall r>0\eqdef{bo}$$
for some $\sigma\in (0,1], \,Q>0$. However, to pass  from inequalities on sets of finite $\nu$ measure to the whole of $\Rn$, we cannot use  the exponential regularization Lemma 9 as is, since we are using two different measures in it. As it will be apparent from the proof below, we need to introduce some conditions at infinity satisfied by the measure $\nu$, in order to regularize the inequality on the whole space. It turns out that it is enough to ask that there are $r_1, Q'>0$ such that 
$$\nu(E)\le Q'|E|,\qquad \forall E {\hbox{ Borel  measurable with } } E\subseteq \{x: |x|\ge r_1\}.\eqdef{cond}$$
This condition is equivalent to asking that, outside a fixed ball, $\nu$ is  absolutely continuous w.r. to the Lebesgue measure, with bounded Radon-Nikodym derivative.
 An example is the singular measure $d\nu(x)=|x|^{(\sigma-1)n}dx$ considered in  [LL1], [LL2], [AY] and other papers.

For simplicity we formulate here only a version of Theorem 1 for these more general measures, however  analogous statements  can be made for other theorems of this paper, for example  Theorem 7,  Theorem  15, or even Theorem 18, in the context of hyperbolic spaces.

\proclaim Theorem 22. Let $\nu$ be a positive Borel measure on $\Rn$ satisfying \eqref{bo}.  If  $0<\alpha<n$ there exists a constant $C=C(\alpha,n,\sigma,Q)$ such that:
 \smallskip\item{(a)} For every measurable and compactly supported  $f:\Rn\to\R$ such that 
$$\|f\|_\na^\na+\|I_\a*f\|_\na^\na\le 1\eqdef{bo1}$$
and for all Borel measurable $E\subseteq \R^n$ with $\nu(E)<\infty$, we have 
$$\int_E \exp\bigg[{\sigma\over|B_1|}|I_\a*f(x)|^\nna\bigg]d\nu(x)\le C(1+\nu(E)).\eqdef{bo2}$$
If in addition $\nu$  satisfies \eqref{cond} then 
$$\int_{\sR^n}\exp_{[\nsa-2]}\bigg[{1\over|B_1|}|I_\a*f(x)|^\nna\bigg]d\nu(x)\le C.\eqdef{bo3}$$
\smallskip\item{(b)} If $\a$ is an   integer  then, for every  $u\in W^{\a,\nsa}(\Rn)$ such that 
$$\|u\|_\na^\na+\|\nabla^\a u\|_\na^\na\le1\eqdef{bo4}$$
and for all Borel measurable $E\subseteq \Rn$ with $\nu(E)<\infty$, we have
$$\int_E\exp\Big[\sigma\gamma_{n,\a}|u(x)|^{n\over n-\a}\Big]d\nu(x)\le C(1+\nu(E)).\eqdef {bo5}$$
If in addition $\nu$ satisfies \eqref{cond} then 
$$\int_{\sR^n}\exp_{[\nsa-2]}\Big[\s\gamma_{n,\a}|u(x)|^{n\over n-\a}\Big]d\nu(x)\le C.\eqdef{bo6}$$
\item{(c)} If there exist $x_0,r_0$ such that 
$\nu( B(x_0,r))\ge c_1 r^{\sigma n}$, for $0<r<r_0$ with $c_1>0$ then the exponential constants in \eqref{bo3}, \eqref{bo6} are sharp. If  there exist $x_0,r_0$ such that 
$\nu(E\cap B(x_0,r))\ge c_1 r^{\sigma n}$, for $0<r<r_0$ with $c_1>0$, then the exponential constants  in \eqref{bo2}, \eqref{bo5} are sharp.
\par

A word of caution:  the measure $\nu$ in this this theorem enters only in the integration of the exponentials. The functions $f,\,I_\a*f,\, u,\,|\nabla^\a u|$ are still in $L^\nsa(\Rn)$ with respect to the Lebesgue measure.\medskip
\pf Proof. The proofs of inequalities \eqref{bo2}, \eqref{bo5} are identical to the corresponding  ones in Theorem 1 for the Lebesgue measure. The point is that the Adams inequality \eqref{12}  in Theorem 7 holds for the measure $\nu$ as above, with the constant $A_g$ replaced by $\sigma^{-1}A_g$. The procedure is exactly the same, except now the main result we use is \eqref{18} of Corollary 6, and the entire proof given in Theorem 7 goes through.

To deal with the \eqref{bo3}, \eqref{bo6} we need to modify Lemma 9. For simplicity we only prove \eqref{bo6} as the other inequality is completely similar. We assume condition \eqref{cond}
and estimate
$$\eqalign{\int_{\sR^n}&\exp_{[\nsa-2]}\Big[\s\gamma_{n,\a}|u(x)|^{n\over n-\a}\Big]d\nu(x)\le\cr& \le\int_{\{u\ge1\}}\exp\Big[\s\gamma_{n,\a}|u(x)|^{n\over n-\a}\Big]d\nu(x)+e^\alpha \int_{\{u\le 1\}}|u(x)|^{\nsa} d\nu(x)\cr}$$
  
now we have 
$$\eqalign{\nu\{u\ge1\}&= \nu\big(\{u\ge 1\}\cap B(0,r_1)\big)+ \nu\big(\{u\ge 1\}\cap B(0,r_1)^c\big)\cr&\le Qr_1^{\sigma n}+Q'|\{u\ge1\}|\le C(1+\|u\|_\na^\na)\cr}$$
so that we can use \eqref{bo5} to estimate the exponential integral over the set $\{u\ge 1\}$.
Finally, 
$$\eqalign{\int_{\{u\le 1\}}|u(x)|^{\nsa}d\nu(x)&\le \nu\big(B(0,r_1)\big)+\int_{\{u\le 1\}\cap\{|x|\ge r_1\}}|u(x)|^\nsa d\nu(x)\cr&\le Qr_1^{\sigma n}+Q'\int_{\{u\le 1\}}|u(x)|^\nsa dx\le C(1+\|u\|_\na^\na)\le C.\cr}$$
The sharpness statements follows as in the proof of the sharpness statement of Theorem 7 (see proof of Corollary 6).\endpf

\bigskip
\centerline{\bf A sharp Trudinger  inequality on bounded domains without boundary conditions}
\bigskip
The next result has to do with smooth and  bounded domains, so it is somewhat unrelated to what we have done so far. We present it here since it is a nice and simple application of Lemma 10.

 Sharp versions of the Trudinger inequality on smooth, connected, bounded domains $\Omega$ for functions $u\in W^{\a,\nsa}(\Omega)$ are only known for $\alpha=1$ [CY], [Ci1], and for $\alpha=2$, if $\Omega$ is a ball [FM2].  In the case $\alpha=1$ Chang-Yang and Cianchi proved that there is $C$ such that
$$\int_\Omega \exp\bigg[2^{-{1\over n-1}}\gamma_{n,1}|u(x)-\overline u|^{n\over n-1}\bigg]dx\le C\qquad u\in W^{1,n}(\Omega),\;\|\nabla u\|_n\le 1,\eqdef{CC1}$$
where $\overline u=|\Omega|^{-1}\int_\Omega u$, and  $\gamma_{n,1}=n\omega_{n-1}^{1\over n-1}$ is the sharp constant for the Moser-Trudinger inequality on $W_0^{1,n}(\Omega)$.

It is clear that some sort of normalization of $u$ is needed, as in \eqref{CC1}, if restrictions are imposed only on  the seminorm  $\|\nabla u\|_n$. Hence, it makes sense to ask about a sharp inequality under the full Sobolev norm condition $\|u\|_n^n+\|\nabla u\|_n^n\le 1$, and with no additional conditions on $u$, in the same spirit as in  the original paper by Trudinger. As far as we know no such result exists, however  we prove here that it can be  easily obtained by combining the Chang-Yang-Cianchi results and Lemma 10.

\proclaim Theorem 23. If $\Omega$ is a smooth, connected and bounded open set in $\Rn$, there exists a constant $C=C(\Omega)$ such that
$$\int_\Omega \exp\bigg[2^{-{1\over n-1}}\gamma_{n,1}|u(x)|^{n\over n-1}\bigg]dx\le C\eqdef{CC2}$$
 for each $u\in W^{1,n}(\Omega)$ with $\|u\|_n^n+\|\nabla u\|_n^n\le 1$. Moreover the exponential constant is sharp.\par
\pf Proof. In Lemma 10 take $\beta'=n$,  $V=Z=\big\{u\in W^{1,n}(\Omega): \ds\int_\Omega  u=0\big\}$, $\;T$ the identity on $V$, and $p(u)=\|\nabla u\|_n$. Then, using that 
$$|\overline u|\le |\Omega|^{-{1\over n}}\|u\|_n=|\Omega|^{-{1\over n}}\big(1-\|\nabla u\|_n^n\big)^{1\over n}$$
we obtain, using Chang-Yang-Cianchi's result together with (B) in Lemma 10
$$\eqalign{&\int_\Omega \exp\bigg[2^{-{1\over n-1}}\gamma_{n,1}|u(x)|^{n\over n-1}\bigg]dx\le\cr&\le\int_\Omega \exp\bigg[2^{-{1\over n-1}}\gamma_{n,1}\Big(|u(x)-\overline u|+|\Omega|^{-{1\over n}}\big(1-\|\nabla u\|_n^n\big)^{1\over n}\Big)^{n\over n-1}\bigg]dx\le Ce^{-\gamma_{n,1}(2|\Omega|)^{-{1\over n-1}}}.\cr}$$
It is not hard to check that  if $0\in\p \Omega$ then the usual Moser sequence 
$$u_\e(x)=\cases{\log{1\over\e} & if $|x|<\e$\cr\log{1\over|x|} & if $\e\le|x|<1$\cr 0 & if  $|x|\ge1$\cr}$$
saturates the exponential constant in  \eqref{CC2}, arguing for example as in [F] pp. 451-453. The point is that as $\e\to0$ we have
 $$\|u_\e\|_{L^n(\Omega)}^n+\|\nabla u_\e\|_{L^n(\Omega)}^n\sim\half \|u_\e\|_{L^n(\sR^n)}^n+\half\| \nabla u_\e\|_{L^n(\sR^n)}^n\sim\half\| \nabla u_\e\|_{L^n(\sR^n)}^n= \half\omega_{n-1} \log{1\over\e}.$$
\endpf

We note that the results in [CY] and [Ci] were also obtained for smooth domains with finitely many conical singularities, in which case the sharp constant is $n(\theta_\Omega)^{1\over n-1}$, where $\theta_\Omega$ is the minimum solid aperture of the cones at the singularities. Needless to say, a result like Theorem 23 also holds under this more general situation, with the same sharp constant $n(\theta_\Omega)^{1\over n-1}$.

\bigskip
\centerline{\scaps 6. Appendix}
\bigskip
In this section we complete the proof of Theorem 3.

Let 
$$L(y)=\bigg(\int_y^\infty \phi(x)^{\b'}dx\bigg)^{\bpi}\le \|\phi\|_{\b'}=\|f^*\|_{\b'}=\|f\|_{\b'}\le1.$$
In what follows we will repeatedly make use of the following inequalities
$$(a+b)^\b\le a^\b+\beta 2^{\b-1}(a^{\b-1} b+b^\b),\qquad ab\le {a^\b\over\b}+{b^{\b'}\over \b'},\qquad a,b\ge0$$
$$\bigg(\sum_1^m a_k\bigg)^\b\le m^\b\sum_1^m a_k^\b,\qquad  a^{\bpi}\le 1+a.$$
Note that if $0\le z_1\le z_2 \le y$ we have
$$\int_{z_1}^{z_2} g(x,y)^\b dx\le\int_{z_1}^{z_2} (1+ H_1(1+|x|)^{-\gamma}\big)^\b dx\le z_2-z_1+ C_3.$$
Also, for $z\ge y$ we have
$$\int_z^\infty g(x,y)^\b dx\le \int_y^\infty e^{{\b\over q}(y-x)} C_2dx={q\over \b}\, C_2=C_4$$
and 
$$\int_{-\infty}^{0} g(x,y)^\b dx\le {\s  }\, J.$$
Next, we note that for $y>0$ H\"older's inequality implies
$$\eqalign{F(y)\ge y-\int_{-\infty}^\infty g(x,y)^\b dx&\ge y-\Big(y+ C_3+C_4+{\s  } J\Big)\cr&\ge- C_3-C_4-{\s  } J.\cr}\eqdef{2}$$
From now on let
$$d^*=C_3+C_4+{\s} J.$$
\def\imp{\Longrightarrow}
Now let for $\lambda\in\R$
$$E_\lambda=\{y\ge 0 : F(y)\le \lambda\}$$
and let us prove that there exists $C_5$ such that 
$$|E_\lambda|\le C_5 (|\lambda|+d^*)\eqdef{3a}$$
Proceeding as in [A] and [FM1], it's enough to prove that there exists $C_6$ such that 
for any $\lambda\in\R$ 
$$y_1,y_2\in E_\lambda,\quad y_2>y_1>C_6(|\lambda|+d^*)\imp y_2-y_1\le C_6(|\lambda|+d^*)\eqdef{3}$$
indeed, if that is the case, then
$$\eqalign{|E_\lambda|&=\big|E_\lambda\cap\{y: y\le C_6(|\lambda|+d^*)\}\big|+\big|E_\lambda\cap\{y: y> C_6(|\lambda|+d^*)\}\big|\cr& \le  C_6(|\lambda|+d^*)+\sup_{y_2>y_1>C_6(|\lambda|+d^*)\atop y_1,y_2\in E_\lambda} (y_2-y_1)\cr&\le 3C_6(|\lambda|+d^*)\cr} $$
which implies {3a}.

We now prove \eqref3. If $y_1,y_2\in E_\lambda$ and $|\lambda|<y_1<y_2$, then $F(y_2)\le\lambda$, so that 

$$\eqalign{(y_2-\lambda)^\bi&\le\int_{-\infty}^\infty g(x,y)\phi(x)dx=\int_{-\infty}^{y_1}+\int_{y_1}^{y_2}+\int_{y_2}^\infty\cr&\le
\bigg(\int_{-\infty}^{y_1}g(x,y)^\b dx\bigg)^{\bi}+
\bigg[\bigg(\int_{y_1}^{y_2} g(x,y)^\b dx\bigg)^{\bi}+\bigg(\int_{y_2}^\infty g(x,y)^\b dx\bigg)^{\bi}\bigg] L(y_1)\cr&\le
(y_1+d^*)^\bi +\big[(y_2-y_1+C_3)^{\bi}+C_4^\bi\big]L(y_1)
\cr}$$
from which we deduce
$$\eqalign{&y_2-\lambda\le y_1+d^*+\b 2^{\b-1}\big[(y_1+d^*)^\bpi\big((y_2-y_1+C_3)^\bi+C_4^\bi\big)L(y_1)+\cr&\hskip22em + \big((y_2-y_1+C_3)^\bi+C_4^\bi\big)^\b L(y_1)^\b\Big]\cr&\le y_1+d^*+\b 2^{\b-1}\big[(y_1+d^*)^\bpi\big((y_2-y_1+C_3)^\bi+C_4^\bi\big) L(y_1)+2^\b (y_2-y_1+C_3)L(y_1)^\b+2^\b C_4 \Big].
\cr}\eqdef 4$$
Now we show that there exists $C_7$ such that 
$$(y_1+d^*)L(y_1)^{\b'}\le C_7\big(|\lambda|+d^*\big)\eqdef 5$$
Indeed, proceeding as above\def\bp{{\b'}}
$$\eqalign{y_1-\lambda&\le\bigg(\int_{-\infty}^\infty g(x,y_1)\phi(x)dx\bigg)^\b=\bigg(\int_{-\infty}^{y_1}+\int_{y_1}^\infty\bigg)^\b\le\Big(
(y_1+d^*)^\bi\big (1-L(y_1)^{\b'}\big)^\bpi+C_4^\bpi L(y_1)\Big)^\b
\cr
&\le (y_1+d^*)\big (1-L(y_1)^{\b'}\big)^{\b\over\bp}+\b 2^{\b-1}\big[(y_1+d^*)^\bpi\big (1-L(y)^{\b'}\big)^{\b-1\over\bp}C_4^\bi L(y_1)+C_4 L(y_1)^\b\big]\cr&
\le(y_1+d^*)\Big (1-{\b\over\bp}L(y_1)^{\b'}\Big)+C_4^\bi\b 2^{\b-1}(y_1+d^*)^\bpi L(y_1)+\b 2^{\b-1}C_4
\cr}$$
or 
$$-\lambda\le d^*-{\b\over\b'}(y_1+d^*)L(y_1)^\bp+C_4^\bi \b2^{\b-1}(y_1+d^*)^\bpi L(y_1)+\b2^{\b-1}C_4.$$
Letting $z=(y_1+d^*)^\bpi L(y_1)$ the last inequality can be written as
$$ z^\bp\le C_8(z+ \lambda+d^*)\le {C_8^\b\over\b}+{ z^{\b'}\over\b'}+C_8(|\lambda|+d^*)$$
which proves \eqref5. Back to \eqref 4 
$$\eqalign{y_2-y_1&\le \lambda+d^*+(y_2-y_1+C_3)^\bi\big[\b 2^{\b-1}(y_1+d^*)^\bpi L(y_1)\big]+\b2^{\b-1}C_4^\bi(y_1+d^*)^\bpi L(y_1)\cr&\hskip17em+\b2^{2\b-1}(y_2-y_1+C_3)L(y_1)^\b+C_4 \b2^{2\b-1}\cr
&\le \lambda+ d^*+{y_2-y_1+C_3\over\b}+{\big(\b2^{\b-1}\big)^\bp(y_1+d^*)L(y_1)^\bp\over\b'}+\b2^{\b-1}C_4^\bi C_7^\bpi(|\lambda|+d^*)^\bpi\cr&\hskip17em +\b2^{2\b-1}(y_2-y_1+C_3)L(y_1)^\b+C_4 \b2^{2\b-1}
\cr&\le {y_2-y_1\over \b}+C_9(|\lambda|+d^*)+C_{10}(y_2-y_1)L(y_1)^\b
\cr
}$$
so that
$${y_2-y_1\over\bp}\le C_9(|\lambda|+d^*)+C_{11}(y_2-y_1)\bigg({|\lambda|+d^*\over y_1+d^*}\bigg)^{\b\over\bp}. $$
Taking $y_1>2\b' C_{11}(|\lambda|+d^*)$ gives 
$y_2-y_1\le 2\b'C_{11}\big(|\lambda|+d^*\big), $ which is \eqref 4.

\smallskip
To complete the proof we now estimate
$$\eqalign{\int_0^\infty &e^{-F(y)}dy=\int_{-d^*}^\infty |E_\lambda|e^{-\lambda}d\lambda\le \int_{-d^*}^\infty \big(C_5(|\lambda|+d^*)\big)e^{-\lambda}d\lambda\le C_{12}d^*e^{d^*}\le C_{13}\Big(1+{\s }J\Big) e^{{\s }J}. 
\cr}$$
\endpf
\smallskip
\centerline{\bf References }
\vskip1em
\item{[AT]} Adachi S., Tanaka K., 
{\sl Trudinger type inequalities in $\R^N$ and their best exponents}, 
 Proc. Amer. Math. Soc. {\bf 128} (2000), 2051-2057.\smallskip 
\item{[A]} Adams D.R. {\sl
A sharp inequality of J. Moser for higher order derivatives},
Ann. of Math. {\bf128} (1988), no. 2, 385--398. 
\smallskip
\item{[AY]} Adimurthi, Yang Y., {\sl  An interpolation of Hardy inequality and Trundinger-Moser inequality in $\R^N$ and its applications}, IMRN {\bf 13} (2010), 2394-2426.\smallskip

\item{[BGS]} Banica V., Gonz\'alez M.d.M., S\'aez M., {\sl Some constructions for the fractional Laplacian on noncompact manifolds}, Rev. Mat. Iber., to appear.\smallskip

\item{[Cao]} Cao D. M., 
{\sl Nontrivial solution of semilinear elliptic equation with critical exponent in $\R^2$}, 
Comm. Partial Differential Equations {\bf  17} (1992),  407-435. \smallskip
\item{[CST]} Cassani D., Sani F., Tarsi C.,
{\sl Equivalent Moser type inequalities in $\R^2$ and the zero mass case}, J. Funct. Anal. {\bf 267} (2014), 4236-4263.\smallskip 
\item{[CY]} Chang S.-Y.A., Yang  P.C., {\sl  Conformal deformation of metrics on $S^2$}, J. Differential Geom. {\bf 27} (1988), 259-296.\smallskip
\item{[Ci1]} Cianchi A., {\sl Moser-Trudinger trace inequalities}, Adv. Math. {\bf 217} (2008), 2005-2044.\smallskip
\item{[Ci2]} Cianchi A., {\sl Moser-Trudinger inequalities without boundary conditions and isoperimetric problems}, Indiana Univ. Math. J. {\bf 54} (2005), 669-705. \smallskip
\item{[CK]} Cohl H.S., Kalnins E.G.,  {\sl Fundamental solution of the Laplacian in the
hyperboloid model of hyperbolic geometry}, preprint (2012) arXiv:1201.4406.\smallskip
\item{[CL]} Cohn W.S., Lu G., {\sl Best constants for Moser-Trudinger inequalities on the Heisenberg group}, Indiana Univ. Math. J. {\bf50} (2001), 1567-1591. \smallskip
\item{[do\'O]}  do \'O J.M.B. {\sl $N$-Laplacian equations in $\R^N$ with critical growth}, Abstr. Appl. Anal. {\bf 2} (1997), 301-315.\smallskip
\item{[F]}  Fontana L.,  {\sl Sharp borderline Sobolev inequalities on compact Riemannian manifolds}, Comment. Math. Helv. 
{\bf68} (1993), 415--454.\smallskip
\item{[FM1]} Fontana L., Morpurgo C., {\sl Adams inequalities on measure spaces}, Adv. Math. {\bf 226} (2011), 5066-5119.
\smallskip
\item{[FM2]} Fontana L., Morpurgo C., {\sl Sharp Moser-Trudinger inequalities for the Laplacian without boundary conditions}, J. Funct. Anal. {\bf 262} (2012),  2231-2271. 
\smallskip
\item{[IMN]} Ibrahim S., Masmoudi N., Nakanishi K., {\sl Trudinger-Moser inequality on the whole plane with the exact growth condition}, JEMS, to appear.\smallskip
\item{[IMM]} Iula S., Maalaoui A., Martinazzi L., {A fractional Moser-Trudinger type inequality in one dimension and its critical points}, preprint (2015), arXiv:1504.04862. \smallskip
\item{[LL1]} Lam N., Lu G., {\sl 
Sharp Moser-Trudinger inequality on the Heisenberg group at the critical case and applications}, Adv. Math. {\bf 231} (2012),  3259-3287.\smallskip 
\item{[LL2]} Lam N., Lu G., {\sl A new approach to sharp Moser-Trudinger and Adams type inequalities: a rearrangement-free argument}, 
J. Differential Equations {\bf 255} (2013), 298-325.\smallskip 
\item{[LLZ]} Lam N., Lu G., Zhang L., {\sl Equivalence of critical and subcritical sharp Trudinger-Moser-Adams inequalities},  arXiv:1504.04858 (2015).\smallskip
\item{[LR]} Li Y., Ruf B., {\sl A sharp Trudinger-Moser type inequality for unbounded domains in $R^n$},   
 Indiana Univ. Math. J. {\bf 57} (2008),  451-480.\smallskip
\item{[LT1]}  Lu G., Tang  H., {\sl  Best constants for Moser-Trudinger inequalities on high dimensional hyperbolic spaces}, Adv. Nonlinear Stud. {\bf 13} (2013),  1035-1052. \smallskip
\item{[LT2]}  Lu G., Tang  H., {\sl  Sharp Moser-Trudinger Inequalities on Hyperbolic
Spaces with Exact Growth Condition}, J. Geom. Anal., to appear.\smallskip
\item{[LTZ]} Lu G., Tang  H., Zhu M., {\sl Best constants for Adams' inequalities with the exact growth condition in $\R^n$}, preprint (2015).\smallskip
 \item{[MS]} Mancini G., Sandeep K., {\sl  Moser-Trudinger inequality on conformal discs}, Commun. Contemp. Math. {\bf 12} (2010),  1055-1068. \smallskip
\item{[MST]} Mancini G., Sandeep K., Tintarev C.,
{\sl Trudinger-Moser inequality in the hyperbolic space $\H^N$}, 
 Adv. Nonlinear Anal. {\bf 2} (2013),  309-324.\smallskip
\item{[MS1]} Masmoudi N.; Sani F.,{ \sl Trudinger-Moser inequalities with the exact growth condition in $\R^n$} and applications, Comm. Partial Differential Equations, to appear.\smallskip
\item{[MS2]}  Masmoudi N., Sani F., {\sl Adams' inequality with the exact growth condition in $\R^4$}, Comm. Pure Appl. Math. {\bf 67} (2014), 1307-1335.\smallskip
\item{[Mo]} Moser J., {\sl A sharp form of an inequality by N. Trudinger}, Indiana Univ. Math. J. {\bf20} (1970/71), 1077-1092.\smallskip

\smallskip
\item{[MST]} Mancini G., Sandeep K., Tintarev C., {\sl Trudinger-Moser inequality in the hyperbolic space $\H^N$},
Adv. Nonlinear Anal. {\bf 2} (2013), 309-324.\smallskip 
\item{[Og]} Ogawa T. {\sl 
A proof of Trudinger's inequality and its application to nonlinear Schrödinger equations}, 
Nonlinear Anal. {\bf 14} (1990),  765-769.\smallskip 

\item{[Oz]}  Ozawa T., {\sl On critical cases of Sobolev's inequalities}, J. Funct. Anal. {\bf 127} (1995), 259-269.\smallskip 
\item{[Pa]} Panda R. {\sl 
Nontrivial solution of a quasilinear elliptic equation with critical growth in $\R^n$}, Proc. Indian Acad. Sci. Math. Sci. {\bf 105} (1995), 425-444.\smallskip 

\item{[Ruf]} Ruf B., {\sl A sharp Trudinger--Moser type inequality for unbounded domains in $\R^2$}, J. Funct. Anal. {\bf 219} (2005), 340-367.\smallskip
\item{[RS]} Ruf B., Sani F., {\sl
 Sharp Adams-type inequalities in $\R^n$},  
Trans. Amer. Math. Soc. {\bf 365} (2013), 645-670.\smallskip 
\item{[S]} Stein E.M., {\sl
Singular integrals and differentiability properties of functions}, 
Princeton Mathematical Series {\bf 30}, Princeton University Press, 1970.\smallskip
\item{[Str1]} Strichartz R.S., {\sl Multipliers on fractional Sobolev spaces}, J.  Math. Mech. {\bf 16} (1967), 1031-1060.\smallskip
\item{[Str2]} Strichartz R.S., {\sl A note on Trudinger's extension of Sobolev's inequalities}, Indiana Univ. Math. J. {\bf 21} (1971/72), 841-842.\smallskip
\item{[Tar]} Tarsi C., 
{\sl Adams' inequality and limiting Sobolev embeddings into Zygmund spaces}, 
Potential Anal. {\bf 37} (2012),  353-385.\smallskip 
\item{[Tat]}  Tataru D., {\sl Strichartz estimates in the hyperbolic space and global existence for the semilinear wave equation}, Trans. Amer. Math. Soc. {\bf 353} (2001), 795-807.\smallskip
\item{[Tr]} Trudinger N.S., {\sl 
On imbeddings into Orlicz spaces and some applications}, 
J. Math. Mech. {\bf17} (1967), 473-483. 
\smallskip
\item{[Y]} Yang Y., {\sl 
Trudinger-Moser inequalities on complete noncompact Riemannian manifolds},  
J. Funct. Anal. {\bf 263} (2012),  1894-1938.\smallskip

\smallskip\smallskip
\noin Luigi Fontana \hskip19em Carlo Morpurgo

\noin Dipartimento di Matematica ed Applicazioni \hskip5.5em Department of Mathematics 
 
\noin Universit\'a di Milano-Bicocca\hskip 12.8em University of Missouri

\noin Via Cozzi, 53 \hskip 19.3em Columbia, Missouri 65211

\noin 20125 Milano - Italy\hskip 16.6em USA 
\smallskip\noin luigi.fontana@unimib.it\hskip 15.3em morpurgoc@missouri.edu

\end